\input amstex
\documentstyle{amsppt}
\magnification1200
\tolerance=1000
\def\n#1{\Bbb #1}
\def\fr{\hbox{fr }}

\def\Tr{\hbox{Tr }}
\def\cl{\hbox{ cl }}

\def\loc{\hbox{loc }}
\def\im{\hbox{im }}

\def\Gal{\hbox{Gal }} 
\def\Hom{\hbox{Hom }}
\def\Ker{\hbox{Ker }}
\def\good{{good}}

\def\Res{\hbox{Res }}
\def\Endom{\hbox{End }}
\def\Sel{\hbox{Sel }}
\def\Inv{\hbox{Inv }}
\def\Div{\hbox{Div }}
\def\Aut{\hbox{Aut }}
\def\Alb{\hbox{Alb }}
\def\Br{\hbox{Br }}

\def\bad{{bad}}

\def\Supp{\hbox{ Supp }}

\def\diag{\hbox{ diag }}
\def\ve{\varepsilon}
\def\vf{\varphi}
\topmatter
\title 
Reduction of a problem of finiteness of Tate-Shafarevich group to a result of 
Zagier type 
\endtitle 
\author
Dmitry Logachev 
\endauthor
\NoRunningHeads 
\date November, 2004 \enddate
\thanks I am grateful to B. Edixhoven, G. Faltings, R. Weissauer, Th. Zink for 
numerous consultations on the subject of the present paper, especially of 
Sections 4, 6. At the inicial stages of preparation of this paper I visited 
universities of Heidelberg, Paris-Sud, Grenoble, Paris-Nord, Bielefeld, 
Rennes, M\"unster, and IHES. I am grateful to L.Clozel, B.Edixhoven, M.Flach, 
E.Freitag, R.Gillard, P.Schneider, J.Tilouine, Th.Zink for their kind 
hospitality. I am grateful to B.H.Gross, S.Kudla, J.Nekov\'a\v r, R.Taylor 
for some important remarks, and to Th.Berry for linguistic corrections. 
\endthanks
\abstract Kolyvagin proved that the Tate-Shafarevich group of an elliptic 
curve over $\n Q$ of analytic rank 0 or 1 is finite, and that its algebraic 
rank is equal to its analytic rank. A program of generalisation of this 
result to the case of some motives which are quotients of cohomology motives 
of Shimura varieties is offered. We prove some steps of this program, mainly 
for quotients of $H^7$ of Siegel sixfolds. For example, we ``almost'' find 
analogs of Kolyvagin's trace and reduction relations. Some steps of the 
present paper are new contribution, because they have no analogs in the case 
of elliptic curves. There are still a number of large gaps in the program. 
The most important of these gaps is a high-dimensional analog of a result of 
Zagier about ratios of Heegner points corresponding to different imaginary 
quadratic fields on a fixed elliptic curve. The author suggests to the readers 
to continue these investigations. \endabstract 
\address 
Departamento de Matem\'aticas 
Universidad Sim\'on Bol\'\i var
Apartado Postal 89000 
Caracas, Venezuela
\endaddress
\endtopmatter
\document 
{\bf 0. Introduction.}
\medskip
Let $E$ be an elliptic curve over $\n Q$ of analytic rank 0 or 1. Kolyvagin 
([K89], [K90] and subsequent papers) proved the
\medskip
{\bf Theorem 0.1.} (a) $SH(\n Q,E)$ --- the Tate-Shafarevich group of $E$ 
over $\n Q$ --- is finite; 

(b) the rank of $E(\n Q)$ is equal to the analytic rank of $E$. 
\medskip
There is the following problem 
\medskip
{\bf 0.2.} Generalise (0.1) to the case of some motives which are quotients 
of cohomology motives of Shimura varieties. 
\medskip
It turns out that (0.2) is a very difficult problem. The present paper is the 
third in a series of papers (the first two papers are [L04], [L05]) whose 
purpose is 
\medskip
(1) To offer a program of a proof of (0.2) (quoted below as The Program); 
\medskip
(2) To prove some steps of The Program, especially for submotives of $H^7$ 
of a Siegel sixfold $X$. The main unconditional result is Theorem 2.13. The 
main conditional result is described in 5.5. 
\medskip
The Program follows the ideas of the original Kolyvagin's proof of (0.1), with some modifications and new steps; for the most essential new contributions see 2.1a,b and Section 3. Since the main property of analogs of Euler systems for the present case is weaker than in [K89], [K90], these analogs are called pseudo-Euler systems. 

It is necessary to emphasize that for 4 different types of Shimura varieties $X$ and Hecke correspondences $\goth T_p$ on $X$: 

(a) $X$ is a Siegel variety of genus 2, $\goth T_p=T_p$; 

(b) $X$ is a Siegel variety of genus 3, $\goth T_p=T_p$; 

(c) $X$ is a Siegel variety of genus 3, $\goth T_p=T_{p,1}$; 

(d) $X$ is a Siegel variety of genus $\ge 4$, $\goth T_p=T_p$

we get obstacles of 4 different types to the realization of The Program (one obstacle for one type of $X$, $\goth T_p$). See Remark 1.6a (and also [L01]) for the type (a), Section 4.4 for the type (b), and obstacles of types 2c, 2d (see below) for the types (c), (d) respectively. Maybe for the latter two cases these obstacles will be got over. 

Are there cases where there is no obstacles? I don't know. 
\medskip
Let us describe now the steps of The Program which are not made yet. They can be subdivided into 2 types. 
\medskip
Problems of type 1 are well-known conjectures of general mathematical interest concerning the main objects of the present paper: 
\medskip
(0.3) Langlands conjecture for Siegel sixfolds $X$;

(0.4) Conjectures about existence and properties of quotient motives of $H^7(X)$; 

(0.5) Serre conjecture on the image of $l$-adic representations in their cohomology;

(0.6) Problem of construction and properties of smooth compact models of $X$. 
\medskip
The origin of the problems of type 2 is The Program itself. They can be subdivided into 4 subtypes. 
\medskip
Problems of type 2a are of purely technical nature. They are easy but time-consuming. Their list is given in the main text. 
\medskip
Problems of type 2b are more complicated, but without doubt solvable. For example, these are problems of rigorous proof of properties of reductions of some subvarieties of Siegel varieties (Sections 4, 6).
\medskip
Problems of type 2c exist thanks to a phenomenon which does not occur in the proof of (0.1), namely the existence of the so-called bad components (see [L05]). A rough analog of the first problem of type 2c in the 1-dimensional case is a calculation of Gross-Kohnen-Zagier ([Z85], [GKZ87]) of ratios of Heegner points corresponding to different imaginary quadratic fields on a fixed elliptic curve over $\n Q$. 
\medskip
In our case analogs of 
\medskip
{\it Heegner points corresponding to different imaginary quadratic fields} 
\medskip
are 
\medskip
{\it some cycles on $X$}. 
\medskip
These cycles depend on a prime $p$. So, the first unsolved problem is to prove the existence of ratios of Abel - Jacobi images of these cycles for different $p$, and to calculate these ratios. The second problem is to prove existence of $p$ such that this ratio does not satisfy a certain congruence. See (1.7)--(1.10) for details. Maybe these problems will be solved by a method similar to [Z85] (see Appendix 3). 

The fact that the problems of type 2c seem to be the most complicated among other unsolved problems for the case when $X$ is a Siegel sixfold, explains the title of the paper. 

If $X$ is a Siegel variety of genus $\ge 4$ then we cannot find reductions of some Shimura subvarieties of $X$ (relevant to our situation) using methods of Sections 4, 6 of the present paper, although this is necessary for the realisation of The Program. This obstacle is of type 2d. See section 4.5 for details.  
\medskip
The structure of [L04], [L05] and the present paper is the following. Properties of quotient motives of $H^7(X)$ are given in [L04]. These results are conditional, we assume the truth of (0.3), (0.4). [L05] contains a generalization of Kolyvagin's trace relations (see (1.3), (1.4)) for $X$ and a Hecke correspondence on $X$. 
\medskip
Section 1 of the present paper contains a survey of Kolyvagin's proof of (0.1) and the corresponding steps of The Program, together with the detailed description of the present paper. In Section 2 we collect together unconditional steps of The Program obtaining a criterion of finiteness of $SH(\n Q,E)$ and $E(\n Q)$ of an abelian variety $E$ over $\n Q$ (Theorem 2.13). The crucial object of the statement of Theorem 2.13 is an operator $U(p)$; proof of its existence is equivalent to the solution of problems of type 2c. 
\medskip
In Section 3 we give a universal method of construction of cycles on Shimura varieties which are homologically equivalent to 0. 
\medskip
Sections 4 -- 6 give ideas of application of Theorem 2.13 to Siegel sixfolds. Practically, their contents is a generalization of Kolyvagin's reduction relation (1.5). The case $\goth T_p=T_p$ (resp. $\goth T_p=T_{p,1}$) is treated in Section 4 (resp. 6). These results are conjectural. The problem of finding the exact answers to the above problems are problems of type 2b. Moreover, we formulate all propositions as if $E$ is an abelian variety. Really, $E$ is a quotient motive of $H^*(X)$. The problem of rewriting of all propositions in terms of cohomology groups of motives is a problem of type 2a. This is made (for another situation) in [N92].
\medskip
Sections 4.3 -- 4.4 contain a ``counterexample'': the case of a Siegel sixfold and a Hecke correspondence $T_p$. This case is interesting, because in spite of the existence of bad components we are able to get an exact value of $U(p)$. Unfortunately, in this case $U(p)$ does not satisfy condition (2.15b) of Theorem 2.13, so we cannot prove finiteness of SH using $T_p$. 
\medskip
This example is included for 2 reasons. Firstly, there is still the possibility of error in arguments (change of a sign would be sufficient!), which can imply a happy end. The second reason (the main one) is the following: maybe our methods might be applicable in other situations. For example, let us consider the functional case where analogs of abelian varieties are Drinfeld modules and shtukas. We know that the functional case is usually easier. Maybe in some cases either 
\medskip
(a) There will be no bad components, or:
\medskip
(b) The phenomenon of Section 4.4 (the pseudo-Euler systems for the case $g=3$, $\goth T_p=T_p$ are identically 0) does not occur, or: 
\medskip
(c) It will be possible to use the method of [Z85] (see Appendix 3) in order to find the Abel-Jacobi image of bad components? 
\medskip
Further, Section 4.5 contains a case of a Siegel variety of genus 4. We get that probably in this case there is no obstacle of type 2c but a new obstacle of type 2d appears. 
\medskip
Section 5 contains a possible example where there are no trivial arguments that $U(p)$ is always ``bad'': the case of a Siegel sixfold and the Hecke correspondence $T_{p,1}$. 
\medskip
{\bf 1. Survey of Kolyvagin's proof, and parallel steps of the present paper.} 
\medskip
For the convenience of the reader, we give here a survey of ideas of the original Kolyvagin's proof for the case of $E$ of analytic rank 0 (all details of secondary importance are omitted). They will be marked by (*). In parallel, we indicate the corresponding steps of the present paper, they will be marked by (**). We do not consider here the case of $E$ of analytic rank 1, because even the case of rank 0 is rather complicated. We use notations of Kolyvagin and we use [G91] for references. 

{\bf (*)} Let $N$ be a level, $X_0(N)=\overline{\Gamma_0(N)\backslash\Cal H}$ the compactification of the modular curve of level $N$, $\vf: \Alb (X_0(N)) \to E$ a Weil map to an elliptic curve $E$ over $\n Q$. Let $l$ and $M=l^n$ be a fixed prime and its power (both $l$ and $M$ are denoted in [G91] by $p$). Further, let $K$ be an imaginary quadratic field (in order to simplify the proofs we consider only the case when $h(K)=1$). We denote by $\Sel (E/\Bbb Q)_M$ the Selmer group. 

Recall the definition of Heegner point. Points $t$ on the open part of $X_0(N)$ are in one-to-one correspondence with the isogenies of elliptic curves $\psi_t: A_t \to A'_t$ such that $\Ker \psi_t=\n Z/N\n Z$. A point $t\in X_0(N)$ is called a Heegner point with respect to $K$ if both $A_t$, $A'_t$ have complex multiplication by the same order of $K$. A Heegner divisor is a Galois orbit of a Heegner point; Heegner divisors are exactly 0-dimensional Shimura subvarieties of $X_0(N)$ in the sense of Deligne ([D71]). 

$\vf$ is defined on divisors of degree 0 on $X_0(N)$. To transform a Heegner divisor of degree $d$ to a divisor of degree 0, we subtract $d$ times the image of the cusp $i\infty$ on $X_0(N)$. Its $\vf$-image is a Heegner point on $E$. For a given $K$ there exists the ``principal'' Heegner divisor $x_1\in \Div(X_0(N))(K)$ (which is one point if $h(K)=1$) and its $\vf$-image --- the ``principal'' Heegner point $y_1\in E(K)$ which are the main objects of [K89] (see also [G91], first page). 
\medskip
The main result of [K89] is the 
\medskip
{\bf Proposition 1.1.} If $y_1$ is not a torsion point and $\Tr_{K/\n Q}(y_1)\in E(\n Q)$ is a torsion point, then there exists $\goth c$ which does not depend on $l$, $M$ such that $\goth c \Sel (E/\Bbb Q)_M =0$. 
\medskip
Finiteness of $SH(E/\Bbb Q)$ and $E(\Bbb Q)$ follows immediately from this proposition. For simplicity, we shall consider in this survey only the case $\Tr_{K/\n Q}(y_1)=0$. 
\medskip
{\bf (**)} In the most general setting of (0.2) an analog of $X_0(N)$ is any ``modular object'' $X$ and an analog of a Heegner divisor on $X_0(N)$ is a subobject $V$ of $X$. For example, in [KL92] we treat the case when $X$ is a quaternionic Shimura curve. 

Particularly, let $X$ be a smooth compact model of a Shimura variety and $V$ a codimension $d$ cycle on $X$ such that 

(a) $V$ is homologically equivalent to 0; 

(b) The support of $V$ is a union of Shimura subvarieties of $X$ and cycles with support at infinity; 

(c) $X$ and $V$ are defined over a number field $k$. 

There exists the $l$-adic Abel - Jacobi image of $V$
$$\cl '(V)\in
H^1(k,H^{2d-1}_{\hbox {et}}(X\otimes \bar{\n Q},\n Z _l(d)))$$ 
Let $E$ be a quotient motive of $H^{2d-1}(X)$; it is an analog of the elliptic curve $E$ of [K89]. We can prolonge $\cl'$ to the $l$-adic cohomology of $E$, this is an analog of $\vf$ of [K89]. 

Definitions of analogs of $SH$ and of the rank of $E(\n Q)$ are given in [BK90]. The main theorem 2.13 of the present paper is formulated for the case when $E$ is an abelian variety. Since quotients of cohomology motives of Shimura varieties which are treated in the present paper are not motives of abelian varieties, we need an analog of Theorem 2.13 for motives. This proof is not given; this problem is of type 2a, it can be solved as in [N92]. Analogously, in Sections 4 -- 6 we treat $E$ as if it were an abelian variety.

Most calculations of the present series of papers are made for the case when $X$ is a smooth compact model of a Siegel sixfold of level $N$, and $V$ is a Picard modular surface (if $h(K)=1$). So, $d=4$ and $E$ is an irreducible quotient motive of $H^7(X)$. 

The definition of $V$ and of the inclusion $V \hookrightarrow X$ is given in [L05]. Recall that points of the open part of $X$ parametrize isogenies of abelian threefolds with kernel $(\n Z/N\n Z)^3$ and points of $V$ correspond to those threefolds whose endomorphism algebra is the maximal order of $K$. 

We do not consider in the present paper problems related to a smooth compact model of $X$. These problems are of type (0.6). 
\medskip
{\bf (*)} We denote $\Cal V=K(E_M)$ the field generated by $M$-torsion points of $E$ (it is denoted by $V$ in [K89] and by $L$ in [G91], Section 9, p. 249). We consider the $l$-adic representation 
$$\rho_l:\Gal(K)\to \Aut(E_M)=GL_{2}(E_M)$$
{\bf (1.1a)} We shall consider in this survey only cases when $\rho_l$ is a surjection. The general case can be easily reduced to this one. 
\medskip
We choose a prime $p$ such that 
\medskip
{\bf (1.2)} The Frobenius of $p$ in $\Cal V/\n Q$ is the complex conjugation.
\medskip
($p$ is denoted in [G91] by $n=l_1\cdot \dots \cdot l_k$ or, if $k=1$, simply by $l$. For the case when the analytic rank of $E$ is 0 we can choose $k=1$, $n=p$ a prime). 

Particularly, $p$ is inert in $K$. Let us recall (for the case $h(K)=1$) the definition of the ring class field $K^p$ of $K$ (denoted by $K_n$ in [G91]): it is the only abelian extension of $K$ with Galois group $\Gal K^p/K =\Bbb Z/(p+1)\Bbb Z$, non-ramified outside $p$, totally ramified at $p$ and such that the corresponding subgroup of the idele group of $K$ contains the idele whose $p$-component is $p$ and other components are 1. We denote $\Gal K^p/K$ by $G=G_p$, and we choose and fix its generator $g=g_p$ (denoted by $\sigma_l$ in [G91]). 

We have the Hecke correspondence $T_p$ on $X_0(N)$. Its restriction to $E$ is multiplication by $a_p$ --- the $p$-th Fourier coefficient of the normalised cusp form of weight 2 corresponding to $E$. 

Attached to $p$ are a Heegner point $x_p\in X_0(N)(K^p)$ and its $\vf$-image --- a Heegner point $y_p\in E(K^p)$ ([G91], p. 238). There are formulas: 
$$\Tr_{K^p/K}(x_p)= T_p(x_1)\eqno{(1.3)}$$
(equality of divisors on $X_0(N)$);
$$\Tr_{K^p/K}(y_p)= a_p y_1\eqno{(1.4)}$$
$$\widetilde{y_p}=\fr (\widetilde{y_1})\eqno{(1.5)}$$
where tilde means reduction at a valuation over $p$ and fr is the Frobenius automorphism of $\bar \n F_p$ ([G91], Proposition 3.7). (1.3), (1.4) are called Kolyvagin's trace relations for $X_0(N)$ and $E$ respectively, and (1.5) is called Kolyvagin's reduction relation. 
\medskip
{\bf (**)} The first step of The Problem is to generalize these trace and reduction relations. The problems of generalization of (1.3), (1.5) are of independent interest regardless of their application to a solution of The Problem for some cases. 

The obtained results are the following. The paper [L05] is devoted to finding of analogs of (1.3) for the case when $X$ is a Siegel variety and $\goth T_p$ a $p$-Hecke correspondence on $X$. There exist 2 finite sets $L_{\good}$, $L_{\bad}$ (depending on $\goth T_p$), and for all $i\in L_{\good}$ (resp. $j\in L_{\bad}$) there are irreducible subvarieties $V_{p,i}$, $V_{p,j}$ defined over $K^p$ (resp. over $K$) such that we have an equality of cycles on $X$: 

$$\goth T_p(V)=\left(\bigcup_{i\in L_{\good}}\alpha_{p,i}\left(\bigcup_{\beta=0}^pg^{\beta}(V_{p,i})\right)\right) \cup \left(\bigcup_{j\in L_{\bad}}\alpha_{p,j}(V_{p,j})\right) \cup \alpha V\eqno{(1.6)}$$
where $\alpha_{p,i}$, $\alpha_{p,j}$, $\alpha$ are multiplicities. [L05] gives the complete answer (i.e. finding of $L_{\good}$, $L_{\bad}$, $V_{p,i}$, $V_{p,j}$, $\alpha_{p,i}$, $\alpha_{p,j}$, $\alpha$) for the case $X$ is a Siegel sixfold, and $\goth T_p=T_p$ is the simplest $p$-Hecke correspondence. Partial answers are obtained for the cases: 

1. $X$ is a Siegel sixfold, $\goth T_p=T_{p,1}$ the $p$-Hecke correspondence defined by the matrix $\diag(1,1,p,p^2, p^2,p)$. 

2. $X$ is a Siegel variety of genus $> 3$, $\goth T_p=T_p$ is the simplest $p$-Hecke correspondence. 
\medskip
{\bf Remark 1.6a.} Inclusion $V \hookrightarrow X$ corresponds to an inclusion of reductive groups $GU(r,s) \hookrightarrow GSp_{2g}$ where $r$, $s$ is the signature of the unitary group, $r+s=g$. It is known that the maximal field of definition of components of $T_p(V)$ is $K^p$ if $r\ne s$ and $K$ if $r=s$, i.e. $L_{\good}\ne\emptyset$ iff $r\ne s$. Existence of good components is a necessary condition for our construction of pseudo-Euler systems. Particularly, for $g=2$ the method of the present paper does not give pseudo-Euler systems. This is why we consider the case $g=3$. 
\medskip
In order to use Abel-Jacobi map, we apply the construction of Section 3 to $V$, $V_{p,i}$, $V_{p,j}$. This construction will give us cycles which are homologically equivalent to 0. Their Abel-Jacobi images are denoted by $y_1$, $y_{p,i}$, $y_{p,j}$ respectively. We treat them as elements of an abelian variety $E$, i.e. $y_1, y_{p,j}\in E(K)$, $y_{p,i}\in E(K^p)$. 
\medskip
The origin of obstacle of type 2c of The Program is the existence of bad components. Since $y_{p,j}\in E(K)$ and for a ``general'' $E$ the rank of $E(K)$ is 1, we can formulate 
\medskip
{\bf Conjecture 1.7.} There exists a coefficient $x_{p,j}\in \n Q$ such that 
$$y_{p,j}=x_{p,j}y_1\eqno{(1.8)}$$
We denote $\goth x_p=\sum_{j\in L_{\bad}}\alpha_{p,j}x_{p,j}$.
\medskip
It turns out that in order to use Theorem 2.13. we must
\medskip
{\bf (1.9).} Find the residue of $x_{p,j}$ (or of $\goth x_p$) modulo $M^2$. 
\medskip
{\bf (1.10).} Prove existence of $p$ (satisfying other conditions of Theorem 2.13) such that $\goth x_p/M$ is not congruent mod $l$ to some number that can be calculated explicitly. 
\medskip
See 5.5 for the final result.
\medskip
{\bf Remark 1.11.} Roughly speaking, we can find $x_{p,j}$ modulo $M$ (Sections 4, 6). For the case $\goth T_p=T_p$ we have: $L_{\bad}$ consists of one element $j_1$, and $\alpha_{p,j_1}=p+1$. Since $p+1$ is a multiple of $M$, knowledge of $x_{p,j}$ modulo $M$ implies knowledge of $x_{p}$ modulo $M^2$. Unfortunately, condition (2.15b) of Theorem 2.13 is not satisfied in this case. 
\medskip
{\bf (*)} Condition (1.2) implies (see for example [G91], (3.3) )
$$M| (p+1)\eqno{(1.12)}$$
$$M| a_p\eqno{(1.13)}$$

Now we consider a commutative square 

$$\matrix E(K)/ M E(K)         & \to & H^1(K,E_M) \\
          \downarrow           &     & \downarrow \\
		  [E(K^p)/ M E(K^p)]^G & \overset{\delta_p}\to{\to} & 
		  [H^1(K^p,E_M)]^G \endmatrix\eqno{(1.14)}$$
(the left square of [G91], (4.2) --- we need only this left square). 
We denote the right vertical map of (1.14) by Res. (1.1a), (1.2) imply that Res is an isomorphism. 

Let $P\in E(K^p)$ be an element such that its image in $E(K^p)/ M E(K^p)$ is $G$-stable. This means that $g(P)-P\in ME(K^p)$. We denote by $c$ the 
element $\Res^{-1}(\delta_p(P))\in H^1(K,E_M)$

(1.1a), (1.2) imply that $E_M \cap E(K^p)=0$. This means that the element 
$$B=\frac{g(P)-P}{M}\in E(K^p)\eqno{(1.15)}$$ 
is well-defined. 

We can identify $E_M$ and $\tilde E_M$. Since $g$ acts trivially on $\tilde E$, we have 
$$\tilde B \in \tilde E_M=E_M\eqno{(1.16)}$$

Let us consider the localization of $c$ at $p$. We denote by $K_p$, $K^p_p$ localizations at $p$ of $K$, $K^p$ respectively. Let $K_p^{(M)}$ be the maximal abelian extension of $K$ such that $\Gal(K_p^{(M)}/K_p)$ is an $M$-torsion group, and let $K_p^{p,(M)}$ be the subfield of $K^p_p$ of degree $M$ over $K_p$. We restrict $g\in \Gal(K^p/K)$ to an element of $\Gal(K_p^{p,(M)}/K_p)$ which we denote by $g$ as well. 

{\bf (1.16a)} Since $K_p^{(M)}/K_p$ is the composite of the disjoint extensions $K_p^{p,(M)}/K_p$ and $\n Q_{p^{2M}}/K_p$ --- the non-ramified extension of degree $M$ of $K_p=\n Q_{p^{2}}$, we can consider $g\in\Gal(K_p^{p,(M)}/K_p)$ as an element of $\Gal (K_p^{(M)}/K_p)$. Further, we denote the Frobenius automorphism of $\n Q_{p^{2M}}$ over $K_p$ by $\fr_{K_p}$, and we can consider $\fr_{K_p}$ as an element of $\Gal (K_p^{(M)}/K_p)$ as well. 
\medskip
Formula (1.2) implies that 

$H^1(K_p,E_M)= \Hom (\Gal (K_p), E_M)= \Hom (\Gal (K_p^{(M)}/K_p), E_M)$. 

We denote by $\loc_p$ the localization map 

$H^1(K,E_M) \to H^1(K_p,E_M)=\Hom (\Gal (K_p^{(M)}/K_p), E_M)$. 

This means that $\loc_p(c)(g)\in E_M$ is well-defined, and we have the following formula of purely cohomological nature (it follows immediately from the reduction of [G91], (4.6)): 

$$\loc_p(c)(g) = \tilde B\eqno{(1.17)}$$

Now we apply the above formulas to the element 
$$D=D_p=\sum_{i=0}^p ig ^i(y_p) \in E(K^p)\eqno{(1.18)}$$ 
([G91], (4.1); notation of [G91] is $P_n$). (1.4), (1.12), (1.13) imply that $D_p \in [E(K^p)/ M E(K^p)]^G$; the corresponding 
$B$, $c$ are denoted by $B_p$, $c(p)$. The element $c(p)\in H^1(K,E_M)$ is an element of level 1 of an Euler system. (1.17) becomes 
$$\loc_p(c(p))(g) = \tilde B_p\eqno{(1.19)}$$
\medskip
Now we consider the image of $y_1$ in $E(K)/ M E(K) \hookrightarrow H^1(K,E_M)$ and denote it by $c(1)$. (1.16a) shows that $\loc_p(c(1))(\fr_{K_p})\in E_M=\tilde E_M$ is well-defined. We have: 
$$\tilde B_p = -\fr(\loc_p(c(1))(\fr_{K_p}))\eqno{(1.20)}$$
where the first $\fr\in \Gal(\bar \n F_p/\n F_p)$ acts on $\tilde E_M$. (1.20) follows from the definitions of $D_p$, $B_p$, $c(1)$, formulas (1.4), (1.5), (1.12), (1.13) and the formula for the characteristic polynomial of Frobenius on $\tilde E$: $\fr^2-a_p \fr +p=0$. See [G91], calculations on the upper half of page 246. So, we have a formula 
$$\loc_p(c(p))(g)=-\fr(\loc_p(c(1))(\fr_{K_p}))\eqno{(1.21)}$$
(the main property of Euler systems of level 1). 
\medskip
{\bf Remark 1.21a.} Since for high-dimensional cases the characteristic polynomial of Frobenius on $\tilde E$ is distinct from $\fr^2-a_p \fr +p=0$, (1.20) and (1.21) do not hold in high-dimensional cases. 
\medskip
Now let $s\in \Sel (E/\Bbb Q)_M\hookrightarrow H^1(K,E_M)$ be any element. We want to show that $s=0$ (some non-essential multipliers that contribute to $\goth c$ of (1.1) are neglected). We consider the Tate pairing 
([G91], 7.3) of $s$ and $c(p)$: $$<s, c(p)>\in \Br(K)$$ The global Tate pairing 
is the sum of local ones. The local Tate pairing of 2 non-ramified elements 
is 0. The sum of invariants of elements of $\Br(K)$ is 0, $c(p)$ is 
non-ramified at all points of $K$ except $p$, and $s$ is non-ramified at 
all points of $K$. This means that the local Tate pairing of $s$ and $c(p)$ 
at $p$ is 0: 
$$<\loc_p(s),\loc_p(c(p))>=0\eqno{(1.22)}$$

There exists a formula for the local Tate pairing: if $s_1, s_2\in H^1(K_p,E_M)=\Hom (\Gal(K_p^{(M)}/K_p), E_M)$ and $s_1$ is non-ramified, then we have (after some 
identification of $\frac1M\Bbb Z/\Bbb Z$ and the group of $M$-th roots of 1 depending on a choice of $g$) 
$$\Inv (<s_1,s_2>)=[s_1(\fr_{K_p}), s_2(g)]\eqno{(1.23)}$$
where $g$ and $\fr_{K_p}$ are from (1.16a), and $[*,*]$ is the Weil 
pairing. See [G91], (7.6). Applying (1.23) to the case $s_1=\loc_p(s)$, $s_2=\loc_p(c(p))$ we get from (1.21), (1.22): 
$$[\loc_p(s)(\fr_{K_p}),\loc_p(c(1))(\fr_{K_p})]=1\eqno{(1.24)}$$

{\bf (**)} We indicate in Section 1 a work-around that permits us to prove (1.24) in high-dimensional cases, in spite of Remark 1.21a. From now on a large segment of proof of (0.1) (formulas 1.24 --- 1.37) coincides with the corresponding steps of the proof of Theorem 2.13, with the following difference: in high-dimensional cases $E$ is an abelian variety of dimension $\goth d$ (instead of an elliptic curve). In order to avoid repeating, we give here some necessary modifications and later (in 1.25a) we continue to give the survey of the proof of (0.1) under the assumption that $E$ is an abelian variety. 

As earlier we denote $\Cal V=K(E_M)$, so $H^1(\Cal V,E_M)=\Hom (\Gal \Cal V, E_M)$. We consider the $l$-adic representation 
$$\rho_l:\Gal(K)\to \Aut(E_M)=GSp_{2\goth d}(\n Z/M)$$
{\bf (1.25)} We shall consider only cases when $\rho_l$ is a surjection (analog of (1.1a) in 1-dimensional case). 

After proving the Serre conjecture for $E$, the reduction of the general case to the condition (1.25) is easy; this is a problem of type 2a. 

If $\rho_l$ is a surjection then the restriction map $H^1(K,E_M)\to H^1(\Cal V,E_M)$ is an inclusion (because for $i=1,2$ \ \ $H^i(GSp_{2\goth d}(\n Z/M),E_M)=0$, the proof for $\goth d=1$ in [G91], (9.1) is valid for any $\goth d$) and $\Gal(\Cal V/ K)=GSp_{2\goth d}(\n Z/M)$. 
\medskip
{\bf (*) (1.25a)} There are maps 
$$E(K) \to E(K)/M \to H^1(K,E_M) \to H^1(\Cal V,E_M)=\Hom (\Gal \Cal V, E_M)$$
For any element $\alpha\in E(K)$ or $\alpha\in H^1(K,E_M)$ we denote by $\alpha_{(1)}$ its image in $\Hom (\Gal \Cal V, E_M)$. Throughout the paper $t$ will mean an element of $E(K)$ or $H^1(K,E_M)$. In both cases $\Ker(t_{(1)})$ is a subgroup of $\Gal \Cal V$. 
\medskip
{\bf (1.26)} We denote by $W(t)$ the extension of $\Cal V$ that corresponds to $\Ker(t_{(1)})$. 
\medskip
We can consider $t_{(1)}$ as an injection from $\Gal (W(t)/\Cal V)$ to $E_M$; we denote this injection by $t_{(2)}$. 

We denote by $\sigma$ the complex conjugation. 
\medskip
{\bf Lemma 1.27.} $W(t)/K$ is a Galois extension. If moreover there exists $\ve_t=\pm 1$ such that 
$$\sigma(t)=\ve_t \cdot t\eqno{(1.28)}$$ 
then $W(t)/\n Q$ is a Galois extension. $\square$
\medskip
The Galois group $\Gal(\Cal V/ K)=GSp_{2\goth d}(\n Z/M)$ acts on $\Gal (W(t)/\Cal V)$. 
\medskip
{\bf Lemma 1.29.} $t_{(2)}$ is a $GSp_{2\goth d}(\n Z/M)$-homomorphism (respectively the above action of $GSp_{2\goth d}(\n Z/M)$ on $\Gal (W(t)/\Cal V)$ and the natural action of  $GSp_{2\goth d}(\n Z/M)$ on $E_M$). $\square$
\medskip
{\bf (1.29a.)} Let $\goth g\in \Gal (W(t)/\Cal V) \subset \Gal (W(t)/\n Q)$ and $p \in \n Z$ a prime such that $\fr_p(W(t)/\n Q)=\sigma \goth g$. 
\medskip
{\bf Lemma 1.30.} Such $p$ satisfies (1.12), (1.13). $\square$
\medskip
Further we shall consider only $p$ satisfying (1.29a) for some $\goth g$. 
\medskip
{\bf Lemma 1.31.} If $t\in H^1(K,E_M)$ is non-ramified at $p$ (particularly, if $t\in E(K)$) then $$t_{(2)}((\sigma \goth g)^2)=\loc_p(t)(\fr_{K_p})\ \ \ \square$$

Let $\Gal(\Cal V)^{(M)}$ be the maximal abelian $M$-torsion quotient group of $\Gal(\Cal V)$. For each subset $C\subset \Hom(\Gal(\Cal V), E_M)=\Hom(\Gal(\Cal V)^{(M)}, E_M)$ we consider (following [G91]) an extension $W(C) \supset \Cal V$ that corresponds to a subgroup 
$$H(C) =\bigcap_{h\in C}\Ker h \subset \Gal(\Cal V)$$
{\bf (1.31a).} Now we return to the above $s$, $y_1$. Let $<*,*>$ denote the linear envelope of elements. We take $C=<s_{(1)}, (y_1)_{(1)}>$, and we denote $W(<s_{(1)}, (y_1)_{(1)}>)$ simply by $W$. For $\goth g\in \Gal (W/\Cal V)$ let $\goth g_s$, $\goth g_{y_1}$ be projections of $\goth g$ on $\Gal (W(s)/\Cal V)$, $\Gal (W(y_1)/\Cal V)$ respectively. 
\medskip
We shall use the following version of 1.29a. For $\goth g\in \Gal (W/\Cal V) \subset \Gal (W/\n Q)$ we shall consider primes $p \in \n Z$ such that $\fr_p(W/\n Q)=\sigma \goth g$. 

According (1.31), we get that (1.24) becomes 

$$[s_{(2)}((\sigma \goth g_s)^2), (y_1)_{(2)}((\sigma \goth g_{y_1})^2)]=1\eqno{(1.32)}$$

Both $s\in H^1(K,E_M)$, $y_1\in E(K)$ satisfy (1.28) with $\ve_s=1$, $\ve_{y_1}=-1$. Let us calculate $t_{(2)}((\sigma \goth g)^2)$ for any $t$ satisfying (1.28) and any $\goth g\in \Gal (W(t)/\Cal V)$. Clearly $\sigma$ acts on both $\Gal (W(t)/\Cal V)$ and $E_M$. We have $(\sigma (t_{(2)}))(\goth g)=\sigma(t_{(2)}(\sigma \goth g \sigma ^{-1}))$, hence 
$$t_{(2)}((\sigma \goth g)^2)=t_{(2)}(\goth g)+\ve_t \sigma(t_{(2)}(\goth g))\eqno{(1.35)}$$
The idea of the end of the proof of (0.1) is the following. We choose $\goth g$ such that for both $t=s$, $t=y_1$ we have 
$$\sigma(t_{(2)}(\goth g))=\ve_t t_{(2)}(\goth g)\eqno{(1.36)}$$ 
In this case (1.35) becomes $t_{(2)}((\sigma \goth g)^2)=2t_{(2)}(\goth g)$ and (1.32) becomes (we consider only case $l\ne 2$) 
$$[s_{(2)}(\goth g_s), (y_1)_{(2)}(\goth g_{y_1})]=1\eqno{(1.37)}$$
{\bf (1.38).} Since $y_1$ is not a torsion point, we get that $(y_1)_{(2)}$ is ``far from 0'' (i.e. if $l^k(y_1)_{(2)}=0$ then $k$ is a large number). Since the Weil pairing is non-degenerate and eigenvalues of $s$ and $y_1$ with respect to $\sigma$ are opposite (1 for $s$ and --1 for $y_1$), (1.37) implies that $s_{(2)}$ is ``close to 0'' (i.e. there exists a small number $k$ such that  $l^ks_{(2)}=0$). Practically this implies $s=0$ (we neglect some 
multipliers that contribute to $\goth c$ of proposition 1.1). 

We do not give here a more detailed exposition of the end of [K89], because in Section 2 we give a more general and simple proof suitable for the case when $E$ is an abelian variety. 
\medskip
\medskip   
{\bf 2. Proof of the unconditional theorem. Pseudo-Euler elements. }
\medskip
{\bf (2.1).} As was indicated above, we cannot directly imitate Kolyvagin's proof in the present case, because:

(a) The proof of (1.20) uses a fact that the characteristic polynomial of
Frobenius on a modular curve $X_0(N)$ is

$$\fr^2-T_p \fr +p$$
which does not hold on general $X$.

(b) $\dim H^{2d-1}(E)=2\goth d$ where $\goth d>1$, so arguments related with orthogonality of $s_{(2)}(\goth g_s)$ and $(y_1)_{(2)}(\goth g_{y_1})$ (see (1.38)) must be changed. (The reader might think that we need $\goth d$ independent Heegner elements in order to use arguments of orthogonality; really we need only one).

Now we formulate a theorem that formalizes the situation. Let $E$ be an abelian variety over $\n Q$ of dimension $\goth d$. Let $l$, $M=l^n$, $K$, $\Cal V=K(E_M)$ be as
earlier ($l\ne 2$). We assume that $E$ satisfies (1.25). We denote $\Cal V_{2n}=K(E_{M^2})$. We fix a
simplectic basis $\goth B$ of $E_{M^2}$ over $\n Z/M^2$, i.e. the matrix of the Weil pairing in this basis is the simplectic matrix $J_{2\goth d}=\left(\matrix 0&E_{\goth d}\\ -E_{\goth d}&0\endmatrix\right)$. Let $p$ be a
prime satisfying the following condition 
\medskip
{\bf (2.2).} The matrix of the action of $\fr_p (\Cal V_{2n}/\n Q)$ on $E_{M^2}$ in the basis $\goth B$ is $$\diag(\underbrace{1+aM,\dots,1+aM}_{\goth d\hbox{times}},\underbrace{-1+bM,\dots,-1+bM}_{\goth d\hbox{ times}})$$
where $a$, $b\in \n Z/M$ are some numbers satisfying $a\not\equiv 0$,
$b\not\equiv 0$, $a\not\equiv b \mod l$. 
\medskip
Further, let $y_1 \in E(K)$, $y_p \in
E(K^p)$ be any elements satisfying:
\medskip
{\bf (2.3)} $\Tr_{K/\n Q}(y_1)=0$;
\medskip
{\bf (2.4)} $y_1$ is not a torsion point, and, moreover, $y_1$ is not a multiple of any other element of $E(K)$; 
\medskip
{\bf (2.5)} $\Tr_{K^p/K}(y_p)=\kappa_p y_1$ in $E(K)$, $\kappa_p$ is an integer
coefficient.
\medskip
Now let us consider $W(y_1)$ as in (1.26), and we impose the following condition on $p$ (a stronger version of (1.29a): 
$$\fr_p (W(y_1)/\n Q)=\sigma \goth g\eqno{(2.6)}$$
where $\goth g$ is an element of $\Gal(W(y_1)/\Cal V)$ of order exactly $M$. Later in (2.14) we introduce a stronger version of this condition and prove (Lemma 2.17) that (2.2) and (2.6) are compatible. 
\medskip
So, we can imitate the Kolyvagin's construcion of an element of Euler system as
follows. We define $D_p\in E(K^p)$ like in (1.18). 
\medskip
{\bf Proposition 2.8.} The image of $D_p$ in $E(K^p)/ M E(K^p)$ is $G$-stable.
\medskip
{\bf Proof.} We must prove that $g(D_p)-D_p \in M E(K^p)$. Since
$$g(D_p)-D_p=\Tr_{K^p/K}(y_p)-(p+1)y_p=\kappa_p y_1-(p+1)y_p\eqno{(2.9)}$$
and (2.2) implies that $M|(p+1)$, it is sufficient to prove that $M|\kappa_p$. (2.2) implies $\tilde E(\n F_{p^2})_{l^{\infty}}=\tilde E_M$ (the index $l^{\infty}$ means the $l^{\infty}$-torsion subgroup or the projection of elements to this subgroup). Further, $\tilde D_p \in \tilde E(\n F_{p^2})$. We consider the projection of the reduction of $g(D_p)-D_p$ to $\tilde E(\n F_{p^2})_{l^{\infty}}$. From one side, it is 0, because $\tilde g$ on $\tilde E(\n F_{p^2})$ is trivial. From another side, (2.9) implies that it is equal to $\kappa_p (\tilde y_1)_{l^{\infty}}$. So, in order to prove the proposition, it is sufficient to prove that $(\tilde y_1)_{l^{\infty}}$ has the order exactly $M$. 

Let $(\tilde y_1)_{(3)}$ mean the projection of $\tilde y_1\in \tilde E(\n F_{p^2})$ in $\tilde E(\n F_{p^2})/M$. Condition $(\tilde y_1)_{l^{\infty}}$ has the order exactly $M$ is equivalent to the condition that the order of $(\tilde y_1)_{(3)}$ in $\tilde E(\n F_{p^2})/M$ is exactly $M$.

Now we untroduce some notations for the lemma 2.10 below. Let $\goth E$ be any abelian variety of dimension $\goth d$ over a finite field $\n F_p$ such
that the matrix of the action of $\fr$ on $\goth E_{M^2}$ is 

$$\diag(\underbrace{1+\goth aM,\dots,1+\goth aM}_{\goth d\hbox{
times}},\underbrace{-1+\goth bM,\dots,-1+\goth bM}_{\goth d\hbox{ times}})$$

We denote $\goth E(\n F_{p^2})^-=\{x\in \goth E(\n F_{p^2})|\fr(x)=-x\}$. We define a map $\beta: \goth E(\n F_{p^2})^- \to \goth E_M$ as follows: for $z \in \goth E(\n F_{p^2})^-$ let $z_{(4)} \in H^1(\n F_{p^2}, \goth E_M)$ be its image under the Kummer map. Since $\Gal(\n F_{p^2})$ acts trivially on $\goth E_M$, $z_{(4)}(\fr^2)\in \goth E_M$ is defined. We let $\beta(z)=z_{(4)}(\fr^2)$. 
\medskip
{\bf Lemma 2.10.} In the above notations $\beta(z)=-2\goth bz$. $\square$
\medskip
So, it is sufficient to prove that $(\tilde y_1)_{(4)}(\fr_{\n F_{p^2}})$ is of order $M$. Using (1.31) it is sufficient to prove that $(y_1)_{(2)}((\sigma \goth g)^2)$ is of order exactly $M$. Since $\goth g$ is of order exactly $M$, we see that the left hand side of (1.31) is of order exactly $M$, hence $(\tilde y_1)_{(4)}(\fr_{\n F_{p^2}})$ as well. $\square$ 
\medskip
{\bf Corollary 2.11.} There exists the only element $B\in E(K^p)$ such that
$MB=g(D_p)-D_p$. $\square$
\medskip
{\bf Corollary 2.12.} $\tilde B \in \tilde E_M$. $\square$
\medskip
Now we can formulate the main theorem. Let $E$, $l$, $M$, $K$, $p$, $y_1$, $y_p$, $D$, $B$, $\tilde B$ be as above, $y_1$, $y_p$ satisfy (2.3) - (2.5), $s\in \Sel (E/\Bbb Q)_M$ any element. Let $W=W(<s_{(1)}, (y_1)_{(1)}>)$ be as in (1.31a). 
\medskip
{\bf Theorem 2.13.} If for any $\goth g\in \Gal (W/\Cal V)$ there exists $p$
satisfying (2.2) and the following conditions (2.14), (2.15): 
\medskip
{\bf (2.14)} $\fr_p (W(s,y_1)/\n Q)=\sigma \goth g$; 
\medskip
{\bf (2.15)} There exists an element $U(p)\in \Endom (\tilde E)$ such that
\medskip
(a) $\tilde B = U(p)(\tilde y_1)$;
\medskip
(b) $U(p)|_{\tilde E_M}$ is an isomorphism of $\tilde E_M$; 
\medskip
(c) $U(p)|_{\tilde E_M}$ is a diagonal operator (in the base $\goth B$ restricted on $\tilde E_M$). 
\medskip
Then $s=0$.
\medskip
{\bf Remark 2.16.} Really, the theorem can be proved for the more realistic analogs of conditions (2.3) -- (2.5) and (2.15c). Namely, $\Tr_{K/\n Q}(y_1)$ can be of torsion, $y_1$ can be a multiple of an element of $E(K)$ (in this case $\kappa_p \in \n Q$), $l$ can be 2, etc. These are obstacles of type 2a. 
\medskip
{\bf Proof.} Some steps of the proof coincide with the corresponding steps of 1-dimensional case. We fix an element $\goth g\in \Gal (W/ \Cal V)$ (later in (2.28) we specify $\goth g$) and we consider $p$ satisfying (2.2), (2.14), (2.15) for this $\goth g$.
\medskip
{\bf Lemma 2.17.} Conditions (2.2), (2.14) are compatible.
\medskip
{\bf Proof of 2.17.} It is sufficient to prove that $\Cal V_{2n}$ and $W$ are linearly disjoint over $\Cal V$, because restrictions of (2.2), (2.14) on $\Gal (\Cal V/\n Q)$ coincide.

There exists a subgroup $C_0$ of $\Hom(\Gal(\Cal V), E_M)$ such that $\Cal V _{2n}=W(C_0)$. Really, for each $x\in E_M$ let $\phi _x \in\Hom(\Gal(\Cal V), E_M)$ be defined by the following cocycle formula: $\phi _x (\delta)=\delta(\frac1M x)-\frac1M x$, where $\delta\in \Gal(\Cal V)$ and $\frac1M x$ is fixed. The map $x \mapsto \phi_x$ is a $GSp_{2\goth d}(\n Z/M)$-homomorphism $\phi$ from $E_M$ to $\Hom(\Gal(\Cal V_{2n}/\Cal V), E_M)$. It is clear that $C_0=\phi(E_M)$. 

It is sufficient to show that $<s_{(1)}, (y_1)_{(1)}>\cap C_0=0$ in $\Hom(\Gal(\Cal V)^{(M)}, E_M)$. If there exists $x \in E_M$ such that $\phi_x \in <s_{(1)}, (y_1)_{(1)}>$ then we can assume that $lx=0$. $\Hom(\Gal(\Cal V)^{(M)}, E_M)$ is a $GSp_{2\goth d}(\n Z /M)$-module. Since both $s_{(1)}, (y_1)_{(1)}$ are $GSp_{2\goth d}(\n Z /M)$-stable, the dimension of the linear envelope of $GSp_{2\goth d}(\n Z /M)(\phi_x)$ is $\le 2$. Since $\phi$ is a $GSp_{2\goth d}(\n Z /M)$-homomorphism, the same dimension is ${2\goth d}$ - a contradiction. $\square$ 
\medskip
Let $c(p)$ be as in Section 1 (see lines between (1.18) and (1.19)). Formula (1.19) holds in the present case. Since (1.22), (1.23) also hold, we get 
$$[\loc_p(s)(\fr_{K_p}),\tilde B]=1$$
(2.15a) implies 
$$[\loc_p(s)(\fr^2), U(p)(\tilde y_1)]=1\eqno{(2.20)}$$
(2.15b, c) and (2.20) imply 
$$[\loc_p(s)(\fr^2), (\tilde y_1)_{l^\infty}]=1\eqno{(2.21)}$$
Applying lemma 2.10 to $(\tilde y_1)_{l^\infty}$ we get 
$$[\loc_p(s)(\fr^2), \loc_p(y_1)(\fr^2)]=1\eqno{(2.22)}$$
--- the analog of (1.24) for the present case. 
\medskip
{\bf Remark. } Since the main property of Euler systems (1.21) is not satisfied in the present case, we call elements $c(p)$ elements of pseudo-Euler system. 
\medskip
In the end of the proof $s$ and $y_1$ enter symmetrically in all formulas, so we change notations (as in [K89]) and denote $s=t_1$, $y_1=t_2$; the index $i$ will be 1 and 2. Both $t_i$ satisfy (1.28) with $\ve_1=1$, $\ve_2=-1$. Let $W(t_i)$, $\goth g_{t_i}$ be as in (1.26), (1.31a) respectively. (2.22) implies $$[(t_1)_{(2)}((\sigma \goth g_{t_1})^2),(t_2)_{(2)}((\sigma \goth g_{t_2})^2) ]=0\eqno{(2.23)}$$
(like (1.24) implies (1.32)). 
\medskip
{\bf Lemma 2.24.} $\exists k_i$ such that $\im (t_i)_{(2)}=l^{k_i}E_M$. 
\medskip
{\bf Proof.} Since $(t_i)_{(2)}: \Gal(W(t_i)/\Cal V) \to E_M$ are $GSp_{2\goth d}(\n Z/M)$-homomorphisms, $\im (t_i)_{(2)}$ are $GSp_{2\goth d}(\n Z/M)$-stable subgroups of $E_M$. But $l^kE_M$ are the only $GSp_{2\goth d}(\n Z/M)$-stable subgroups in $E_M$. $\square$
\medskip 
{\bf Lemma 2.25.} $k_2=0$. 
\medskip
{\bf Proof.} $(y_1)_{(2)}=(t_2)_{(2)}$ is of order $M/l^{k_2}$. Since the composite map $\n Z\cdot t_2/M \hookrightarrow E(K)/M \to \Hom (\Gal (W(t_2)/\Cal V), E_M)$ is well-defined and is an inclusion, we get that the image of $y_1$ in $E(K)/M$ is also of order $M/l^{k_2}$. (2.4) implies $k_2=0$. $\square$ 
\medskip
{\bf Lemma 2.26.} $W(t_1)/\Cal V$, $W(t_2)/\Cal V$ are linearly disjoint 
extensions. 
\medskip
{\bf Proof.} If not then $W(t_1) \cap W(t_2) \ne \Cal V$. Let $h$ be a non-trivial element of $\Gal(W(t_1) \cap W(t_2) /\Cal V)$. We have $\sigma (h)=h=-h$: a 
contradiction, because $ l \ne 2$. $\square$ 
\medskip
Let us denote elements of $\goth B$ restricted to $E_M$ by $e_1, \dots, e_{2\goth d}$. Since the matrix of the Weil pairing on $e_1, \dots, e_{2\goth d}$ is $J_{2\goth d}$, we have $[e_1,e_{\goth d+1}]=\zeta_M$ a primitive $M$-th root of 1. \medskip
{\bf Corollary 2.27.} $\exists h \in \Gal(W/\Cal V)$ such that 
$t_1(h)=l^{k_1}e_1$, $t_2(h)=e_{\goth d+1}$. $\square$
\medskip
{\bf 2.28. End of the proof of 2.13.} We take $\goth g$ of the statement of 2.13 equal exactly to this $h$. This $\goth g$ satisfy (1.36) for both $t_i$. Taking into consideration (1.35) and formulas $\sigma(e_1)=e_1$, $\sigma(e_{\goth d+1})=-e_{\goth d+1}$, (2.23) becomes $[l^{k_1}e_1, e_{\goth d+1}]=1$ and hence $k_1=0$, i.e. $(t_1)_{(2)}=0$. Since $E(K)/M \to \Hom (\Gal (W(t_2)/\Cal V), E_M)$ is an inclusion, we get $s=t_1=0$. $\square$ 
\medskip 
The end of the present paper is devoted to an attempt of a construction of $U(p)$. 
\medskip 
\medskip 
{\bf 3. A universal construction of cycles on Shimura varieties which are homologically equivalent to 0.}
\medskip 
Let $k$ be a number field, $X$ a Shimura variety over $k$, $CH^d(X \otimes k)$ the group of codimension $d$ cycles on $X$ modulo rational equivalence defined over $k$ and $CH^d(X \otimes k)_0$ its subgroup of cycles homologically equivalent to 0. Let $E$ be an irreducible quotient motive of $H^{2d-1}_{\hbox {et}}(X\otimes \bar{\n Q},\n Z _l(d))$. For any cycle $\goth V_0 \in CH^d(X \otimes k)_0$ the Abel-Jacobi image cl$_E'(\goth V_0)$ of $\goth V_0$ in $E$ is defined. 

Let now $\goth V \in CH^d(X \otimes k)$ be a cycle. We can associate to $\goth V$ its Abel-Jacobi image in $E$ canonically up to a multiplier using the following construction. 

We denote $r= \hbox{rank}(CH^d(X\otimes k)/CH^d(X\otimes k)_0)$. Let $m$ be a fixed (sufficiently large) prime, $T_m$ the simplest $m$-Hecke operator on $X$ (see (4.0) below), $a_m$ the eigenvalue of $T_m$ on $E$ and 
$Q_m(Z)=\sum_{j=0}^r b_{m,j}Z^j$ the characteristic polynomial of 
the action of $T_m$ on $CH^d(X\otimes k)/CH^d(X\otimes k)_0$, where $Z$ is an independent variable. We denote 
$$\phi _m(\goth V)\overset{\hbox{def}}\to{=}\sum_{j=0}^r b_{m,j}T_m^j(\goth V)$$
Then $\phi _m(\goth V)\in CH^d(X \otimes k)_0$, and its Abel-Jacobi image in $E$ is defined. 
\medskip
{\bf Proposition 3.1.} For different $m$ cl$_E'(\phi _m(\goth V))$ 
are proportional. 
\medskip
{\bf Proof.} We calculate the double sum in 2 different orders: 

$$\hbox{cl}'_E(\sum_{j_1,j_2=0}^r b_{m_1,j_1}b_{m_2,j_2}T^{j_1}_{m_1}
T^{j_2}_{m_2}(\goth V))=
(\sum_{j=0}^r b_{m_1,j}a_{m_1}^j)\hbox{cl}'_E(\phi _{m_2}(\goth V));$$

$$\hbox{cl}'_E(\sum_{j_1,j_2=0}^r b_{m_1,j_1}b_{m_2,j_2}T^{j_1}_{m_1}
T^{j_2}_{m_2}(\goth V))=
(\sum_{j=0}^r b_{m_2,j}a_{m_2}^j)\hbox{cl}'_E(\phi _{m_1}(\goth V)) \ \ \ \ \ \square$$

{\bf Remark 3.2.} The similar construction was used in [KL92], case of $X$ is a quaternion Shimura curve, $z\in X$ a Heegner point. For this case we have $r=1$, $Q_m(Z)=Z-(m+1)$ and $\phi_m(z)$ is the image of $T_m(z)-(m+1)z$ in an irreducible quotient of $\Alb(X)$. This example shows that the
construction above is reasonable. Moreover, since orders of growth 
of $a_m$ and $m+1$ are different we see that $a_m-(m+1) \to \infty$ as $m\to \infty$. Conjecturally, this is true for all cases: 
\medskip
{\bf (3.3)} The order of growth of $b_{m,j}$ and $a_m$ is such that 
$\sum_{j=0}^r b_{m,j}a_m^j$ tends to infinity. 
\medskip
If so then we have the following elementary 
\medskip
{\bf Lemma 3.4.} Let $\goth V_0=\sum _{i \in I} c_i\goth V_i \in CH^d(X)_0$ be a 
linear combination of codimension $d$ Shimura subvarieties of $X$ such that cl$_E'(\goth V_0)\ne 0$. 
Then 
$\exists i \in I$ such that cl$_E'(\phi _m(\goth V_i))\ne 0$. 
\medskip
{\bf Proof.} $\sum _{i \in I} c_i \cl_E'(\phi _m(\goth V_i))=
(\sum_{j=0}^r b_{m,j}a_m^j)\cl_E'(\goth V_0)$.  $ \square $
\medskip
This means that the $\phi _m$-construction of cycles that are 
homologically 
equivalent to 0 is not worse than any other one. 
\medskip
We fix $m$ and we apply this construction to the case $\goth V=V_*$ where $V_*$ are from (1.6), * is some index. We denote the elements cl$_E'(\phi _m(V_*))$ by $y_*$, and we call them the Abel-Jacobi images of $V_*$. The same construction will be applied for the reduced objects $\tilde X$, $\tilde E$ (reduction at $p$). 
\medskip
\medskip 
{\bf Section 4. Counterexample: case of Hecke correspondence $T_p$.} \nopagebreak
\medskip
{\bf 4.0. Definitions. } \nopagebreak
\medskip
The algebra of $p$-Hecke correspondences on a Siegel variety $X$ of any genus $g$ is the ring of polynomials with $g$ generators denoted by $T_p$, $T_{p,1}, \dots, T_{p,g-1}$. They are double cosets corresponding to the diagonal matrices 
\medskip
$$\diag(\underbrace{1,\dots,1}_{g\hbox{ times}}, \underbrace{p,\dots,p}_{g\hbox{ times}})$$ for $T_p$ and 
$$\diag(\underbrace{1,\dots,1}_{g-i\hbox{ times}}, \underbrace{p,\dots,p}_{i\hbox{ times}}, \underbrace{p^2,\dots,p^2}_{g-i\hbox{ times}}, \underbrace{p,\dots,p}_{i\hbox{ times}})$$ for $T_{p,i}$. 

I shall be interested mainly by the case $g=3$, correspondences $T_p$, $T_{p,1}$. The corresponding matrices are $\diag (1,1,1,p,p,p)$, $\diag (1,1,p,p^2,p^2,p)$ respectively. 
\medskip
{\bf Remark.} For $g=3$ and for the Hecke correspondence $T_{p,2}$ the set $L_{\good}$ is empty ([L05]), so we cannot get pseudo-Euler systems using methods of the present paper. 
\medskip
Let $t \in X$ and $A_t$ the corresponding abelian $g$-fold. The set $T_{p}(t)$ is in 1 -- 1 correspondence with the set of maximal isotropic subspaces $W \subset (A_t)_p=(\n F_p)^{2g}$, and the set $T_{p,i}(t)$ is in 1 -- 1 correspondence with the set of isotropic subspaces $W \subset (\n Z/p^2)^{2g}$ such that $W$ is isomorphic to $(\n Z/p^2)^{g-i}\oplus\n F_p^{2i}$. We refer to these subgroups as $W$ of type $T_p$, $T_{p,i}$, and we denote the set of such $W$ by $S_{g}$, $S_{g,i}$ respectively. Finally, we denote $b(n)=\frac{n(n+1)}2$, $G(j,g)(\n F_p)$ the Grassmann variety of $j$-spaces in the $g$-space over $\n F_p$, dimensions are affine, and $G(j,k,g)(\n Z/p^2)$ a generalized Grassmann variety of submodules of $(\n Z/p^2)^{g}$ which are isomorphic to $(\n Z/p)^{k-j} \oplus (\n Z/p^2)^j$ as abstract modules. 
\medskip
Recall that we consider mainly the case $g=3$, $X$ is a Siegel sixfold and $V\subset X$ a Picard modular surface. Some results of Section 3 hold for a more general case $g$ is any number, $V\subset X$ is a subvariety of dimension $g-1$ whose points parametrize abelian $g$-folds having multiplication by the ring of integers of an imaginary quadratic field $K$. 
\medskip
\medskip
{\bf 4.1. Case of ordinary points. }
\medskip
{\bf 4.1.1. Case of one point.} There are correspondences $\Phi_i$ on $\tilde X$ (see (4.1.3) for a definition) such that 
$$\tilde T_p=\sum_{j=0}^g \Phi_j\eqno{(4.1.2)}$$
$\Phi_0$ is the Verschibung 
correspondence and $\Phi_g$ is the Frobenius map. 
Let us fix notations related to the definition of $\Phi_j$. Let
$t \in X(\bar \n Q)$ be an element such that $\tilde A_t$ is
ordinary. Then there exists a fixed isotropic $g$-dimensional subspace $D_g \subset 
(A_t)_p$ enjoying the following property: 
$\tilde \alpha_t:\tilde A_t \to \tilde A'_t$ is the Frobenius map of 
$\tilde A_t$ iff $\Ker \alpha_t = D_g$. Let $t'$ be an element of $T_p(t)$ and $W \subset (A_t)_p$ the corresponding isotropic subspace. We have: 
$$\tilde t' \in \Phi_j(\tilde t) \iff \dim _{\n F_p} W \cap D_g=j \ \ (j=0,\dots,g)\eqno{(4.1.3)}$$
Moreover, for $t'_1, t'_2 \in T_p(t)$ and corresponding subspaces $W_1$, $W_2$ we have: 
$$\tilde t'_1=\tilde t'_2 \iff W_1 \cap D_g=W_2 \cap D_g\eqno{(4.1.4)}$$
{\it Projection. }We denote by $S_{g}(j)$ the set of $W$ such that $\dim _{\n F_p} W \cap D_{g}=j$. Sets $S_{g}(j)$ form a partition of $S_{g}$. There is the natural projection $\pi_j: S_{g}(j) \to G(j, g)$ ($\pi_j(W)=W \cap D_{g}$). The quantity of poins in the fiber of $\pi_j$ is 
$$p^{b(g-j)}\eqno{(4.1.5)}$$
\medskip
{\bf 4.1.6. Case of subvarieties.} Let $\goth V$ be a subvariety of $X$ such that for a generic point $t\in \goth V(K)$ the reduction at $p$ of the corresponding abelian variety $A_t$ is ordinary. There are schemes $\Phi_j(\tilde \goth V)$; we denote their closed subschemes by $\overline{\Phi_j(\tilde \goth V)}$. The Abel-Jacobi images of $\goth V$, $\tilde \goth V$, $\overline{\Phi_j(\tilde \goth V)}$ will be denoted by $\goth y$, $\tilde \goth y$, $\goth y_j$ respectively (recall that we fix $m$ and we use $\phi_m$-construction of Section 3). Clearly reduction commutes with Abel-Jacobi map. $\Phi_j$ act on $\tilde E$. (4.1.5) implies 
$$\Phi_j (\tilde \goth y) = p^{b(g-j)} \goth y_j\eqno{(4.1.7)}$$
and hence 
$$T_p(\tilde \goth y) = a_p\tilde \goth y=\sum_{j=0}^g p^{b(g-j)} \goth y_j\eqno{(4.1.8)}$$ where $a_p$ is the eigenvalue of $T_p$ on $E$. 
\medskip
{\bf Remark.} Considering only the Abel-Jacobi image we loose many information on schemes $\widetilde{\goth T_p(V)}$ and their irreducible components; we take into consideration only the closed support and the depth of these schemes. 
\medskip
{\bf 4.2. Case of non-ordinary points.}
\medskip
We return to our $V\subset X$. The Abel-Jacobi image of $V$ is denoted by $y_1$. For an odd $g$ a generic point $t\in V$ has the property: the reduction of the abelian $g$-fold $A_t$ has the degree of supersingularity 1, i.e. $\# \Supp (\tilde A_t)_p=p^{g-1}$. Particularly, for $g=3$ the $p$-rank of $A_t$ is $2$. 

Recall that all considerations below are conjectural. They are only a first approach to the subject. The rigorous description claims use of another technique. 
\medskip
$\overline{\Phi_j(\tilde V)}$ are reducible: there are subvarieties $\Psi_j(V)\subset \tilde X$ ($j=0,\dots,g-1)$ such that 
$$\overline{\Phi_j(\tilde V)}= \Psi_{j-1}(V) \cup \Psi_j(V) \eqno{(4.2.1)}$$ 
$j=0,\dots,g$, $\Psi_{-1}=\Psi_{g}=\emptyset$. 

We denote the Abel-Jacobi image of $\Psi_j(V)$ by $z_j$. (4.2.1) implies  
$$\Phi_j(\tilde y_1)=p^{b(g-j)+g-j}z_{j-1}+p^{b(g-j)}z_j\eqno{(4.2.2)}$$
($j=0, \dots, g$, $z_{-1}, z_g=0$). 

Clearly that a correct method to find coefficients $p^{b(g-j)+g-j}$, $p^{b(g-j)}$ of (4.2.2) is to calculate dimensions of the corresponding schemes. I did not do it, and the evidence that these coefficients are correct, comes from (4.2.7) below. 

The geometric description of the partition (4.2.1) is the following. We denote by $D_{g-1}^{\bot}$ a $g+1$-dimensional subspace of $(A_t)_p$ which is the kernel of the reduction map $(A_t)_p\to \widetilde{(A_t)}_p$, and by $D_{g-1}$ its dual space. We have $D_{g-1} \subset D_{g-1}^{\bot}$. 
\medskip
{\bf Remark.} $(A_t)_p$ is an $\n F_{p^2}$-space, because $A_t$ has multiplication by $K$ (recall that $p$ is inert in $K$). It is clear that $D_{g-1}$, $D_{g-1}^{\bot}$ are also $\n F_{p^2}$-spaces. 
\medskip
I think that the following analogs of (4.1.3), (4.1.4) hold ($t', t'_1, t'_2, W, W_1, W_2$ are the same, the meaning of $\Psi_j(t)$ is clear; $\Psi_j(V)=\cup_{t \in V} \Psi_j(t)$ ): 

$$\tilde t' \in \Psi_j(t) \iff \dim _{\n F_p} W \cap D_{g-1}=j \ \ (j=0,\dots,g-1)\eqno{(4.2.3)}$$
$$\tilde t'_1=\tilde t'_2 \iff W_1 \cap D_{g-1}=W_2 \cap D_{g-1}\eqno{(4.2.4)}$$
Particularly, both $\Psi_0(t)$, $\Psi_{g-1}(t)$ consist of one point. 
\medskip
{\bf 4.2.5. Partition (*).} Since sets $\Psi_j(t)$ and formula (4.2.3) are conjectural, we change notations and {\it define } $S^{*}_{g}(j)$ as the set of $W$ such that $\dim _{\n F_p} W \cap D_{g-1}=j$, so conjecturally $\Psi_j(t)=S^{*}_{g}(j)$. Sets $S^{*}_{g}(j)$ form a partition of $S_{g}$. There is the natural projection $\pi^*_j: S^{*}_{g}(j) \to G(j, g-1)$ ($\pi^*_j(W)=W \cap D_{g-1}$). The quantity of poins in the fiber of $\pi^*_j$ is 
$$p^{b(g-j)}+p^{b(g-j)-1}\eqno{(4.2.6)}$$
\medskip
This gives us a formula 
$$\cl '_{\tilde E} (\phi_m(\widetilde{T_p(V)}))=a_p\tilde y_1=T_p(\tilde y_1)=\sum_{j=0}^g (p^{b(g-j)}+p^{b(g-j)-1})z_j\eqno{(4.2.7)}$$
Together with (4.1.7), (4.1.8) this gives us coefficients of (4.2.2).
\medskip
{\bf 4.3. Application to irreducible components of $T_p(V)$. }
\medskip
Firstly we develop some ``general theory''. We can unify irreducible good (resp. bad) components $V_{p,i}$ (resp. $V_{p,j}$) from (1.6) if they have equal multiplicities $\alpha_{p,i}$ (resp. $\alpha_{p,j}$). For example, if $\alpha_{p,i_1}=\alpha_{p,i_2}=\dots =\alpha_{p,i_k}$ then we can set $V_{p,\bar i}= V_{p,i_1}\cup V_{p,i_2}\cup \dots \cup V_{p,i_k}$, and analogously for $V_{p,j}$. 

So, we get sets $\bar L_{\good}$, $\bar L_{\bad}$ which are quotient sets of $L_{\good}$, $L_{\bad}$ respectively, and (1.6) becomes 
$$\goth T_p(V)=\left(\bigcup_{\bar i\in \bar L_{\good}}\alpha_{p,\bar i}\left(\bigcup_{\beta=0}^pg^{\beta}(V_{p,\bar i})\right)\right) \cup \left(\bigcup_{\bar j\in \bar L_{\bad}}\alpha_{p,\bar j}(V_{p,\bar j})\right) \cup \alpha V\eqno{(4.3.0)}$$ 
where $V_{p,\bar i}$, $V_{p,\bar j}$ are not necessarily irreducible, but $\alpha_{p,\bar i}$ and $\alpha_{p,\bar j}$ are different. 
\medskip
{\bf Remark.} It is possible a more ``strong'' unification (practically, we can unify all good (resp. bad) components of $\goth T_{p}(V)$. In this case the formula (4.3.7) will be weaken. See Remark 5.6 for a possible application. 
\medskip
Further, we can consider the double union in (4.3.0) as a simple union: 
$$\goth T_{p}(V)= \bigcup_{l\in L}\alpha_lV_{l} \eqno{(4.3.1)}$$
where $L=\bar L_{\good}\times \Gal(K^p/K)\cup \bar L_{\bad}\cup$ \{ the only element corresponding to $V$ in (4.3.0)\}, $V_l$ is one of the sets $g^{\beta}(V_{p,\bar i})$, or $V_{p,\bar j}$, or $V$ itself. (4.3.1) comes from the corresponding decomposition of $S_{g}$: 
$$S_{g}= \bigcup _{l\in L} S'_{g}(l)\eqno{(4.3.2)}$$
namely: for $t \in V$ \ \ \ $W \in S'_{g}(l) \iff $ the corresponding point of $T_{p}(t) \in V_{l} $. The union is disjoint, i.e. sets $S'_{g}(l)$ form a partition of $S_{g}$. 

We denote the Abel-Jacobi image of $V_l$ by $\goth z_l$. Since reduction commutes with Abel-Jacobi map, the Abel-Jacobi image of $\tilde V_l$ is $\tilde \goth z_l$. 
\medskip
{\bf Problem. } What is a formula for $\tilde \goth z_l$? 
\medskip
{\bf Conjecture 4.3.3.} In some cases there exist coefficients $c_{jl} \in \n Q$ such that 
$$\tilde \goth z_l= \sum _{j=0}^g c_{jl}\Phi_j (\tilde y_1);\eqno{(4.3.4)}$$

In order to find $c_{jl}$ we must intersect partitions $(')$ and $(*)$ defined in (4.2.5) and (4.3.2) respectively. 

We fix $l \in L$, and for any $j=0, \dots, g-1$ we consider the set $S^*_{g}(j) \cap S'_{g}(l)$ and the restriction of $\pi^*_j: S^{*}_{g}(j) \to G(j, g-1)$ on it. We denote this restriction by ${\pi'}^*_{j,l}: S^*_{g}(j) \cap S'_{g}(l) \to G(j, g-1)$. Let us restrict ourselves by the case when the following condition holds: 
\medskip
{\bf Condition 4.3.5. } For any $t\in G(j, g-1)$ the quantity of points in the fiber $({\pi'}^*_{j,l})^{-1}(t)$ is the same (does not depend on $t$). 
\medskip
We denote this quantity by $n_{jl}$. 
\medskip
{\bf Conjecture 4.3.6.} We have a formula
$$\tilde \goth z_l=\sum_{j=0}^{g-1} \frac{n_{jl}}{\alpha_l} z_j\eqno{(4.3.7)}$$
Using (4.2.2) we get easily $c_{jl}$. 
\medskip
Now let us apply the above theory to our case $g=3$. We have ([L05]) $L_{\good}$ consists of 2 elements $i_k$, ($k=1,2$), $\alpha_{p,i_k}=1$, $L_{\bad}$ consists of 1 element\footnotemark \footnotetext{We use here the gothic $j$ in order to avoid confusion with the index of $\Phi_j$.} $\goth j_1$, $\alpha_{p,\goth j_1}=p+1$, and $\alpha=0$. So, $\bar L_{\good}$ consists of 1 element $\bar i_1$, and $\bar L_{\bad}=L_{\bad}$. Identifying $\Gal(K^p/K)$ with $\{0, \dots, p\}$ we get 
$$L=\{0, \dots, p\} \cup \{\goth j_1\}$$
For $0\in L$ we denote $V_0=V_{p,\bar i_1}$ simply by $V_p=V_{p,\good}$ and its Abel-Jacobi image $\goth z_0$ by $y_p=y_{p,\good}$; analogously $\goth z_{\goth j_1}$ is denoted by $y_{p,\bad}$. 

[L05] contains the description of $S'_{g}(l)$ used in the proof of the following proposition. 

{\bf Proposition 4.3.8.} Condition (4.3.5) holds for this case, and numbers $n_{jl}$ are given in the following table:
\medskip
\settabs 3 \columns
\+ & $n_{jl}$, $l\in\{0, \dots, p\}$ & $n_{j,\goth j_1}$ \cr
\medskip
\+ $j=0$ & $p^5-p^3$ & $p^4+p^3$ \cr
\medskip
\+ $j=1$ & $p^2$ & 0 \cr
\medskip
\+ $j=2$ & 0 & $p+1$ \cr
\medskip
{\bf Corollary 4.3.9.} According (4.3.7) we get formulas: 
$$\tilde y_{p} = (p^5-p^3)z_0+p^2z_1\eqno{(4.3.10)}$$
$$\tilde y_{p, \bad} = p^3z_0+z_2\eqno{(4.3.11)}$$

{\bf Proof of 4.3.8.} Recall ([L05]) that the partition (4.3.2) for $g=3$ is the following: for $\goth j_1\in L$ $$W\in S_3'(\goth j_1) \iff \dim _{\n F_{p^2}} \n F_{p^2}W=2$$ and $$W\in \bigcup_{i\in\{0, \dots, p\}\subset L}S_3'(i)\iff \dim _{\n F_{p^2}} \n F_{p^2}W=3$$ Since $g=3$, the space $D_{g-1}$ of (4.2) is $D_2$. We need a lemma: 
\medskip
{\bf Lemma 4.3.12.} We have: $\Supp(\tilde V_{bad})=S_3^*(0)(V) \cup S_3^*(2)(V)$; $\Supp(\tilde V_{p})=S_3^*(0)(V) \cup S_3^*(1)(V)$. 
\medskip
{\bf Proof.} If $W \supset D_{2}$ then $W \subset D_{2}^{\bot}$, so $\n F_{p^2} W = D_{2}^{\bot}$, $\dim _{\n F_{p^2}} \n F_{p^2}W=2$ and hence $W\in S_3'(\goth j_1)$, i.e. $W$ corresponds to $V_{bad}$. 
Inversely, let us consider $W$ corresponding to $V_{bad}$. This means that $\dim _{\n F_{p^2}} \n F_{p^2}W=2$. $D_2 \cap \n F_{p^2}W$ is an $\n F_{p^2}$-space, it can have dimension 0 or 1. If $\dim D_2 \cap \n F_{p^2}W=0$ then $D_2 \cap W=0$ and $t \in S_3^*(0)$. If $\dim D_2 \cap \n F_{p^2}W=1$, i.e. $D_2 \subset \n F_{p^2}W$, then it is easy to see that $D_2 \subset W$. Really, $\n F_{p^2}W$ is an $\n F_{p^2}$-space of dimension 2 which contains its orthogonal. Such spaces contain only one isotropic $\n F_{p^2}$-space of dimension 1, namely their orthogonal. This means that $D_2=(\n F_{p^2}W)^{\bot}$, i.e. $\n F_{p^2}W=D_2^{\bot}$, $W \subset D_2^{\bot}$ and hence $D_2 \subset W$. $\square$
\medskip
Now we calculate the quantities of spaces $W$ of each type. There are $p+1$ elements in $S_3^*(2)$ (all $W$ such that $D_2 \subset W \subset D_2^{\bot}$, they form a $P^1(\n F_p)$). There are $p^4+2p^3+p^2$ elements in $S_3^*(1)$. Really, there are $p+1$ possible lines $W \cap D_2$. We fix such a line. There are $(p^2+1)(p+1)$ isotropic planes in $(W \cap D_2)^{\bot}/(W \cap D_2)$, each of them gives us a $W$. It is necessary to subtract $p+1$ planes that contain $D_2/(W \cap D_2)$, we get $p^3+p^2$ planes and multiply this number by $p+1$. 

Since the total number of points in the bad part (i.e. in $S_3'(\goth j_1)$ ) is $p^4+p^3+p+1$, the quantity of points in the bad part of $S_3^*(0)$ is $p^4+p^3$ and in the good part of $S_3^*(0)$ is $p^6+p^5-p^4-p^3$. $\square$
\medskip
{\bf 4.4. Finding of $U(p)$. }
\medskip
Now we apply Theorem 2.13 to this situation. (4.3.0) becomes 
$$T_p(V)=\bigcup_{g\in \Gal(K^p/K)}g(V_{p,\good})\bigcup (p+1)V_{p,\bad}$$ Taking Abel-Jacobi image we get 
$$a_p y_1 = \Tr_{K^p/K}(y_p) + (p+1)y_{p,\bad}\eqno{(4.4.1)}$$ 
where $a_p$ is the eigenvalue of of $T_p$ on $E$. We take $D_p=\sum_{i=0}^p ig ^i(y_p)$, $B_p = \frac{g(D_p)-D_p}M$ from Section 2. Since 
$$g(D_p)-D_p= (p+1)y_p-\Tr_{K^p/K}(y_p)\eqno{(4.4.2)}$$ we have 
$$g(D_p)-D_p= (p+1)y_p - a_p y_1 + (p+1)y_{p,\bad}\eqno{(4.4.3)}$$ We shall see in (4.4.7) that (2.2) implies $M|a_p$, $M|(p+1)$, hence (according (2.5))
$$\tilde B_p = \frac{p+1}M \tilde y_p - \frac{a_p}M \tilde y_1 +\frac{p+1}M \tilde y_{p,\bad}= \frac{p+1}M \tilde y_p - \frac{a_p}M \tilde y_1 +\frac{p+1}M \kappa_p\tilde y_1\eqno{(4.4.4)}$$
Formulas (4.2.2), (4.3.10), (4.3.11) permit us to represent $\tilde y_p$, $\tilde y_{p,\bad}$ as linear combinations of $\Phi_i(\tilde y_1)$. We can easily find the action of $\Phi_i$ on $\tilde E$ using formulas of [L04] (remark: maybe $\Phi_i$ of [L04] is $\Phi_{g-i}$ of the present paper). 

Let us recall the notations (some letter are made gothic in order to avoid confusion with notations of the present paper). Let $\goth T$ be the subgroup of diagonal matrices in $\goth G=GSp_{2g}$ and $\goth M \subset \goth G$ be the subgroup whose $g\times g$-block structure is $\left(\matrix A&0\\ 0&(A^t)^{-1}\endmatrix\right)$, so $\goth T \subset \goth M \subset \goth G$. There are Satake inclusions of Hecke algebras (see, for example, [FCh], Chapter 7 for general formulas or [L04] for explicit formulas): $$\n H(\goth G) \overset{S_G}\to{\hookrightarrow} \n H(\goth M) \overset{S_T}\to{\hookrightarrow} \n H(\goth T) \hookrightarrow \n Z[U_i^{\pm 1}, V_i^{\pm 1}]$$ where $U_i$, $V_i$ ($i=1,\dots,g$) are independent variables. $\n H(\goth G)(\n Z_p)$ is the algebra of $p$-Hecke correspondences on $X$ and $\n H(\goth M)(\n Z_p)$ is the algebra of $p$-Hecke correspondences on $\tilde X$. Particularly, $\Phi_i\in \n H(\goth M)$. 

Let $\Cal M$ be a ``generic'' irreducible submotive of $X$ of middle weight (i.e. the weight of $\Cal M$ is $b(g)$ --- the dimension of $X$) and $\goth E$ its field of coefficients. We can identify a basis over $\goth E \times \n Q_l$ of the cohomology space $H^{b(g)}(\Cal M)$ with the set of subsets of $1,\dots,g$. For $I \subset \{1,\dots,g\}$ we denote by $\goth f_I$ the corresponding element of this basis and we denote $U_I=\prod_{i\in I}U_i\prod_{i\not\in I}V_i\in \n H(\goth T)$. 
We have: $$S_T(\Phi_i)=\sum_{\#(I)=i}U_I\eqno{(4.4.5)}$$ and 
$S_T\circ S_g(T_p)=\sum_{i=0}^g \Phi_i=\sum_{I\subset 1,\dots,g}U_I=\prod_{i=1}^g(U_i+V_i)$, where $T_p\in \n H(\goth G)$. 

The action of $\n H(\goth M)$ on $H^{b(g)}(\Cal M)$ comes from the following action of $\n Z[U_i^{\pm 1}, V_i^{\pm 1}]$ on $H^{b(g)}(\Cal M)$: 

$$\matrix U_i(\goth f_I)= a_0^{1/g} \goth f_I & \hbox { if } & i \in I \\ 
U_i(\goth f_I)= a_0^{1/g}b_i \goth f_I & \hbox { if } &  i \not\in I \\ 
V_i(\goth f_I)= a_0^{1/g}b_i \goth f_I & \hbox { if } &  i \in I \\ 
V_i(\goth f_I)= a_0^{1/g} \goth f_I  & \hbox { if } & i \not\in I \endmatrix\eqno{(4.4.6)}$$
where $a_0$, $b_i$ ($i=1,\dots,g$) are Weil numbers. They satisfy $a_0^2\prod_{i=1}^g b_i=p^{b(g)}$. 

We are interested in the case of a ``generic'' irreducible submotive $\Cal M^-$ of $X$ of weight $b(g)-1$. A basis of $H^{b(g)-1}(\Cal M^-)$ can be identified with the set of $\goth f_I$ such that $I \subset \{2,\dots,g\}$, formulas of the action of $\n Z[U_i^{\pm 1}, V_i^{\pm 1}]$ on $H^{b(g)-1}(\Cal M^-)$ are as above, and $b_1=p$ (see, for example, [L04], Theorem 4.3). 

For our case $g=3$, a basis of $H^{b(g)-1}(\Cal M^-)$ is $\goth f_{\emptyset}$, $\goth f_2$, $\goth f_3$, $\goth f_{23}$ where $\emptyset$, 2, 3, 23 are subsets of (2, 3). These vectors are the eigenvectors of Frobenius with eigenvalues $a_0$, $a_0b_2$, $a_0b_3$, $a_0b_2b_3$ respectively. $T_p$ acts by multiplication by $a_p=a_0(p+1)(b_2+1)(b_3+1)$. Comparing these eigenvalues with the ones from (2.2) we get 
$$a_0\equiv 1+aM, \ \ b_2 \equiv 1, \ \ b_3 \equiv -1 +(a+b)M \mod M^2, \ \ b_1=p \equiv -1 \mod M\eqno{(4.4.7)}$$
(really, we need these congruences only modulo $M$). This means that $a_p \equiv 0 \mod M^2$, hence the term $\frac{a_p}M \tilde y_1$ of (4.4.4) is 0. 

To find $\tilde y_p$, $\tilde y_{p,\bad}$ we use (4.3.10), (4.3.11) and (4.2.2), (4.4.5) --- (4.4.7). Unfortunately we get $\tilde B_p=0$. Pseudo-Euler system does not exist in this case. Calculations are given in Appendix 1. 
\medskip
{\bf 4.5. Case $g=4$.} 
\medskip
If $g$ is even then for a generic point $t\in V(K)$ the abelian variety $\tilde A_t$ is ordinary. Formulas for this case are similar to the formulas of sections 4.2, 4.3, but easier. The space $D_g\subset (A_t)_p$ is an $\n F_{p^2}$-space of dimension $\frac g2$. Analog of the partition (*) of (4.2.5) is the following. We denote by $S^{*}_{g}(j)$ the set of $W$ such that $\dim _{\n F_p} W \cap D_{g}=j$. Sets $S^{*}_{g}(j)$ form a partition of $S_{g}$. There is the natural projection $\pi^*_j: S^{*}_{g}(j) \to G(j, g)$ ($\pi^*_j(W)=W \cap D_{g}$). Further, for any $l \in L$ and for any $j=0, \dots, g$ we consider the set $S^*_{g}(j) \cap S'_{g}(l)$ and the restriction of $\pi^*_j: S^{*}_{g}(j) \to G(j, g)$ on it. We denote this restriction by ${\pi'}^*_{j,l}: S^*_{g}(j) \cap S'_{g}(l) \to G(j, g)$. 

For $g=4$ we have ([L05], idea of the proof): both $\bar L_{\good}$ and $\bar L_{\bad}$ consist of one element which are denoted by $\bar i_1$ and $\bar j_1$ respectively, $\alpha_{p, \bar i_1}=1$, $\alpha_{p, \bar j_1}$ and $\alpha$ are not calculated yet ($\alpha$ is the quantity of isotropic $\n F_{p^2}$-planes in $(A_t)_p=(\n F_{p^2})^4$). So, (4.3.0) becomes (notations of 4.3):
$$T_{p}(V)=\left(\bigcup_{\beta=0}^pg^{\beta}(V_{p,\good})\right)\cup \alpha_{p,\bar j}V_{p,\bad}\cup \alpha V\eqno{(4.5.1)}$$ 
By analogy with the cases $g=2$ ([L01]) and $g=3$, $\goth T_p=T_p$ we can expect that $\alpha_{p, \bar j_1}$ is a multiple of $p+1$. If this is really true then --- as in section 4.4 --- we do not need to know $\goth x_p$ modulo $M^2$ but only $\goth x_p$ modulo $M$. 

Unfortunately even for $g=4$ (as well as for $g>4$) the analog of condition 4.3.5 does not hold. Partition (4.3.2) is ``almost known'' for any $g$. Really, (4.3.2) is a subpartition of the following partition $$S_g=\bigcup_{i=g/2}^gS_g''(i)\eqno{(4.5.2)}$$ where $S_g''(i)$ is the set of maximal isotropic $W\subset (\n F_{p})^{2g}$ such that $\dim_{\n F_{p^2}}(\n F_{p^2}W)=i$ ([L05]). We have $$S_g''(g)=\bigcup_{l\in \bar L_{\good}\times \Gal(K^p/K)}S_g'(l)\eqno{(4.5.3)}$$ If $W\in S_g''(g)$ then $\dim_{\n F_{p}} W\cap D_g$ can be $0,1,\dots, \frac g2$, and $\dim_{\n F_{p}} \n F_{p^2}(W\cap D_g)=2\dim_{\n F_{p}} W\cap D_g$. Particularly, this means that for $g=4$ \ \ \ $\dim_{\n F_{p}} W\cap D_4$ can be $0,1,2$. If $\dim_{\n F_{p}} W\cap D_4=1$ then $W\cap D_4$ can be any $\n F_p$-subspace of $D_4$ of dimension 1. But if $\dim_{\n F_{p}} W\cap D_4=2$, then $W\cap D_4$ is not an $\n F_{p^2}$-space (for $W\in S_g''(g)$). This means that if $t\in G(2,D_4)$ corresponds to an $\n F_{p^2}$-subspace of $D_4$, then the fiber $({\pi'}^*_{g,l})^{-1}(t)$ is empty for $l$ corresponding to $S_g''(g)$, and the same fiber is non-empty if $t\in G(2,4)$ corresponds to an $\n F_{p}$-plane in $D_4$ which is not an $\n F_{p^2}$-line --- 4.3.5 fails. 

What to do? We see that the support of $\Phi_2(\tilde V)$ is a reducible variety: its 2 irreducible components correspond to cases whether $W\cap D_4$ is or is not an $\n F_{p^2}$-subspace of $D_4$. Can we find formulas of Abel-Jacobi images of these components in terms of $\tilde y_1$ and the action of $\Phi_i$ on it? Maybe it is possible to find these Abel-Jacobi images by explicit calculations, because our varieties are over finite fields? 
\medskip
\medskip
{\bf Section 5. Correspondence $T_{p,1}$ --- a possible example.}
\medskip
This section is an analog of Subsection 4.4. For the case $\goth T_p=T_{p,1}$ \ \ $L_{\good}$, $L_{\bad}$, $\alpha_{p,i}$, $\alpha_{p,j}$ and the partition 
$$S_{3,1}= \bigcup _{l\in L} S'_{3,1}(l)$$
(analog of (4.3.2)) are not known completely. We have ([L05]): 

$L_{\good}$ is not empty, it consists of elements $i_1\dots, i_{k_{g,p}}$. $k_{g,p}$ and the numbers $\alpha_{p,i_1},\dots,\alpha_{p,i_{k_{g,p}}}$ are unknown (their finding is reduced to a large but easy calculation). 
\medskip
{\bf Conjecture 5.0.} All $\alpha_{p,i_k}$ are equal. 
\medskip
This conjecture is suggested by [L05], (5.3.1). We shall assume it; see Remark 5.6 if it is wrong. So, (5.0) implies that $\bar L_{good}$ consists of one element. We denote it by $\bar i$ and the common value of $\alpha_{p,i_k}$ by $\alpha_{p,\bar i}$.  
\medskip
$L_{\bad}$ consists of at least 2 elements $j_0$, $j_1,\dots, j_{k_{b,p}}$ ($k_{b,p}$ is unknown; conjecturally, for large $p$ \ \ $k_{b,p}=1$). $\alpha_{p,j_0}=1$, other $\alpha_{p,j_i}$ are unknown. $\alpha=p^5+p^2$ (maybe $p^4+p$). 

Sets $S'_{3,1}(l)$ are known for $l=j_0$, for $l$=\{ the only element corresponding to $V$ in (1.6)\}, and for $l\in \bar L_{\good} \times \Gal(K^p/K)$ sets $S'_{3,1}(l)$ are known up to a choice of one of two possibilities (see [L05], Conjecture 4.2.20). 

Let $y_p=y_{p,\good}=\sum_{\gamma=0}^{k_{g,p}}y_{p,i_{\gamma}}$, $y_{p,\bad}=\sum_{\gamma=0}^{k_{b,p}} \alpha_{p,j_{\gamma}} y_{p,j_{\gamma}}$. Assuming Conjecture 1.7 we can write $y_{p,\bad}=\goth x_p y_1$. Let $D_p$, $B_p$ be as in (4.4). Analog of (4.4.1) for the present case is ($a_{p,1}$ is the eigenvalue of $T_{p,1}$ on $E$)
$$a_{p,1} y_1 = \alpha_{p,\bar i}\Tr_{K^p/K}(y_p) + \goth x_{p} y_1 + (p^5+p^2) y_1 \eqno{(5.1)}$$ i.e. $$\alpha_{p,\bar i}\Tr_{K^p/K}(y_p)=\kappa_p'y_1\eqno{(5.2)}$$ where $\kappa_p'=a_{p,1}-\goth x_{p}-(p^5+p^2)$. 

We can expect that $\alpha_{p,\bar i}=\goth R(p)$ where $\goth R(X)$ is an unknown polynomial. (2.2) implies $p\equiv -1 \mod M$; using this congruence we can find a number $\eta$ (depending only on $\goth R(X)$) such that $\alpha_{p,\bar i}=M^{\eta}\cdot\gamma$, $(\gamma, M)=1$. The same arguments as in the proof of (2.8) show us that $M|\kappa_p'$, so if $\eta>0$ then we can divide (5.2) by $M$ and repeat the process $\eta$ times getting $$\gamma\Tr_{K^p/K}(y_p)=\kappa_py_1\eqno{(5.3)}$$ where $\kappa_p=\kappa_p'/M^{\eta}$. Since $(\gamma, M)=1$, this is practically (2.5). 

We can conjecture that in this situation we can use Theorem 2.13: there is no ``trivial'' obstacles like in the case of $T_p$, because at least one of $\alpha_{p,j_*}$ --- namely, $\alpha_{p,j_0}=1$ --- is coprime to $M$. Formula (2.2) gives us the residue of $a_{p,1} \mod M^2$ (this can be done easily using formulas of [L04]; particularly, $a_{p,1} \equiv 1 \mod M$). So we get that $\goth x_{p} \equiv 1 \mod M$. 
\medskip
{\bf Remark.} There is an independent method to find $\goth x_{p} \mod M$ using (6.2.4), (6.3.3), (6.3.5) for $i=1$, $l=j_0, \dots, j_{k_{b,p}}$. Comparing this method with the result $\goth x_{p} \equiv 1 \mod M$ we can check formulas of Section 6.
\medskip
The analog of (4.4.4) is 
$$\tilde B_p= \frac{p+1}{M}\tilde y_p - \frac{\kappa_p}{\gamma M}\tilde y_1 \eqno{(5.4)}$$

In order to find $\tilde y_p$ we need to find coefficients $c_{jk\ l_0}$ in representation $\tilde y_p= \sum _{j,k} c_{jk\ l_0} \Phi_j\Phi_k (\tilde y_1)$ (formula 6.3.3 for good components $l_0$). This is a problem of type 2b for $g=3$ and of types 2c, 2d for $g>3$. Section 6 contains ideas of solution of this problem. 
\medskip
{\bf 5.5.} Now we can summarize our results. Let us assume that in future we shall be able
\medskip
(a) To find for any $p$ the number $\goth x_p$ modulo $M^{\eta+2}$ --- a problem of type 2c. 
\medskip
In this case (5.4) and other formulas of this section imply the existence of $U(p)$ satisfying 2.15a,c (i.e. $\tilde B_p=U(p)(\tilde y_1)$). So, we must only  
\medskip
(b) To prove that for any $\goth g$ of Theorem 2.13 there exists $p$ such that 
$U(p)$ satisfies 2.15b, i.e. $\goth x_p/M^{\eta+1}$ does not satisfy a certain congruence modulo $l$. After finding of $\goth R(p)$ and $c_{jk\ l_0}$ this congruence can be easily written down explicitly. 
\medskip
{\bf Remark 5.6.} If numbers $\alpha_{p,i_k}$ are different then we can take the maximal value of $\eta$ such that $M^{\eta}$ divides all $\alpha_{p,i_k}$, and define $y_{p,good}$ as $\sum_{\gamma=1}^{k_{g,p}} \frac{\alpha_{p,i_{\gamma}}}{M^{\eta}} y_{p,i_{\gamma}}$. In this case also there is no trivial reasons for $U(p)$ to be  a non-isomorphism on $\tilde E_M$. 
\medskip
{\bf 6. Idea of finding of $\tilde y_p$.} 
\medskip
{\bf Remark 6.0.} We consider here the case of the Hecke correspondence $\goth T_p=T_{p,i}$. The level of na\"{\i}vit\'e of all considerations of this section is higher than the one of section 4. Moreover, (6.2.9) shows that some affirmations are definitely false. I do not know how to correct them, and I shall be grateful to anybody who will help me. 
\medskip
{\bf 6.1. Case of ordinary points.} 
\medskip
{\bf Case of one point.} 
\medskip
Analog of (4.1.2) for $T_{p,i}$ ($i>0$) is the following: 

$$\tilde T_{p,i}=\sum_{j\ge0, k-j \ge i}^g R_{k-j}(i)\cdot p^{-b(k-j)}\Phi_j \Phi_k\eqno{(6.1.1)}$$
where $R_g(i) = R_g(i,p)$ is the quantity of symmetric 
$g \times g$-matrices with entries in $\n F_p$ of corank exactly 
$i$. 

Particularly, for $g=3$ 
$$\tilde T_{p,1} = {1\over p}(\Phi_0 \Phi_1+ 
\Phi_1 \Phi_2+ \Phi_2 \Phi_3)+ {p^2 -1\over p^3}(\Phi_0 \Phi_2 + 
\Phi_1 \Phi_3) + {p^3 -1\over p^4}\Phi_0 \Phi_3 \eqno{(6.1.2)}$$

{\bf Question 6.1.3. }What is an analog of (4.1.3) for the present case? 
\medskip
{\bf Attempt of answer.} Like above let $t'$ be an element of $T_{p,i}(t)$, $W\subset (A_t)_{p^2}$ the corresponding isotropic subspace (see (4.0), $D_g \subset (A_t)_{p^2}$ the kernel of Frobenius, $D_g$ is isomorphic to $(\n Z/p^2)^g$ as an abstract module. Roughly speaking, 
$$W \hbox{ corresponds to }\Phi_j\Phi_k \iff W \cap D_g= (\n Z/p)^{k-j} \oplus (\n Z/p^2)^j\eqno{(6.1.4)}$$
equality of abstract modules. (The author has a more detailed description of this situation). 
\medskip
{\it Projection. }We denote by $S_{g,i}(j,k)$ the set of $W$ of type $T_{p,i}$ such that $W \cap D_{g}=(\n Z/p)^{k-j} \oplus (\n Z/p^2)^j$ as abstract modules. Sets $S_{g,i}(j,k)$ form a partition of $S_{g,i}$. There is the natural projection $\pi_{i;j,k}: S_{g,i}(j,k) \to G(j,k,g)(\n Z/p^2)$ ($\pi_{i;j,k}(W)=W \cap D_{g}$). 
\medskip
{\bf Case of subvarieties.} 
\medskip
Let $\goth V$, $\goth y$, $\tilde \goth y$ be as in (4.1.6), i.e. $\goth V$ a subvariety of $X$ such that for a generic point $t\in \goth V(K)$ the reduction at $p$ of the corresponding abelian variety $A_t$ is ordinary, and $\goth y$, $\tilde \goth y$ the Abel-Jacobi images of $\goth V$, $\tilde \goth V$ respectively. We denote the Abel-Jacobi image of $\overline{\Phi_j\circ\Phi_k(\tilde V)}$ by $\goth y_{j,k}$. 
\medskip
{\bf Conjecture 6.1.5.} The analog of (4.1.7) is the following: 
$$\Phi_j\circ\Phi_k(\tilde \goth y) = p^{b(g-j)+b(g-k)} \goth y_{j,k}$$
Substituting (6.1.5) to (6.1.1) we get 
$$T_{p,i}(\tilde \goth y)=a_{p,i}\tilde \goth y=\sum_{j,k\ge0, j+i \le k}^g R_{k-j}(i)\cdot p^{-b(k-j)+b(g-j)+b(g-k)} \goth y_{j,k}\eqno{(6.1.6)}$$

{\bf Conjecture 6.1.7.} For all $g,i,j,k$ the coefficient $R_{k-j}(i)\cdot 
p^{-b(k-j)+b(g-j)+b(g-k)}$ of (6.1.6) is the quantity of poins in the fiber of $\pi_{i;j,k}$. 
\medskip
{\bf Remark. } This conjecture is checked by explicit calculation for the case $i=1$, $g=3$. The explicit formula for this case is 
$$T_{p,1}(\tilde \goth y)=p^8\goth y_{01}+p^3\goth y_{12}+\goth y_{23}+(p^6-p^4)\goth y_{02}+(p^2-1)\goth y_{13}+(p^3-1)\goth y_{03}\eqno{(6.1.8)}$$
\medskip
{\bf 6.2. Case of non-ordinary points. } \nopagebreak
\medskip
Here $V\subset X$ is from (4.0). I can only guess what is the analog of formulas (4.2.1), (4.2.2). Most likely there are subvarieties $\Psi_{j,k}(V)\subset \tilde X$, $j,k= 0, \dots, g-1$ such that 
$$\overline{\Phi_j\circ\Phi_k(V)}=  \Psi_{j,k}(V) \cup \Psi_{j-1,k}(V) \cup \Psi_{j,k-1}(V) \cup \Psi_{j-1,k-1}(V)\eqno{(6.2.1)}$$
{\bf Question 6.2.2. }Is $\Psi_{j,k}(V)=\Psi_{k,j}(V)$? 
\medskip
I think that yes.
\medskip
Like above we denote the Abel-Jacobi image of $\Psi_{j,k}(V)$ by $z_{j,k}$. 
\medskip
{\bf Question 6.2.3. }What are coefficients in the analog (4.2.2)? 
\medskip
We can expect (taking into consideration (6.1.5)) that 
$$\Phi_j\circ\Phi_k(\tilde y_1) = p^{b(g-j)+b(g-k)}[z_{j,k}+p^{g-j}z_{j-1,k}+p^{g-k}z_{j,k-1}+p^{(g-j)+(g-k)}z_{j-1,k-1}]\eqno{(6.2.4)}$$
We denote by {\bf (6.2.5)} the result of substitution of (6.2.4) to (6.1.6). Particularly, for the case $i=1$, $g=3$ we have 
$$\tilde T_{p,1}(\tilde y_1)=p^{10}z_{00}+(p^8+p^7+p^6-p^5)z_{01}+p^4z_{11}+$$
$$+(p^6+2p^5-2p^2)z_{02}+(p^3+p^2+p-1)z_{12}+z_{22}\eqno{(6.2.6)}$$
\medskip
We define $D_{g-1}^{\bot}$, $D_{g-1}$ like in (4.2). In our case $D_{g-1}=(\n Z/p^2)^{g-1}$ as an abstract module. 

{\bf Conjecture 6.2.6a.} The following analogs of (4.2.3), (4.2.4) hold (by analogy with (6.1.4); $t', t'_1, t'_2, W, W_1, W_2$ have the analogous meaning): 
$$\tilde t' \in \Psi_{j,k}(t) \iff W \cap D_{g-1}= (\n Z/p)^{k-j} \oplus (\n Z/p^2)^j\ \ \eqno{(6.2.7)}$$
($k\ge j$; $j,k=0,\dots,g-1$; equality of abstract modules); 
$$\tilde t'_1=\tilde t'_2 \iff W_1 \cap D_{g-1}=W_2 \cap D_{g-1}\eqno{(6.2.8)}$$
\medskip
{\bf Partition (*).} We denote by $S^{*}_{g,i}(j,k)$ the set of $W$ of type $T_{p,i}$ such that $W \cap D_{g-1}=(\n Z/p)^{k-j} \oplus (\n Z/p^2)^j$ as abstract modules. Sets $S^{*}_{g,i}(j,k)$ form a partition of $S_{g,i}$. There is the natural projection $\pi^{*}_{i;j,k}: S^{*}_{g,i}(j,k) \to G(j,k,g-1)(\n Z/p^2)$ ($\pi^{*}_{i;j,k}(W)=W \cap D_{g-1}$). 
\medskip
{\bf 6.2.9.} It is natural to expect that the coefficient at $z_{j,k}$ in (6.2.5) for any $g,i,j,k$ is equal to the quantity of poins in the fiber of $\pi^{*}_{i;j,k}$, but this is not true. Explicit calculation of the quantity of poins in the fiber of $\pi^{*}_{i;j,k}$ for the case $i=1$, $g=3$ is given in Appendix 2. This calculation shows that this is true for all pairs $(j,k)$ except the pair (0,2): the quantity of poins in the fiber of $\pi^{*}_{i;j,k}$ is (see theorem A2.4, $(j, \nu)=(0,2)$, where $\nu=k-j$, and $\mu=0,1,2$)
$$p^6+2p^5-p^3-2p^2\eqno{(6.2.10)}$$
while the coefficient at $z_{02}$ is (see 6.2.6)
$$p^6+2p^5-2p^2\eqno{(6.2.11)}$$
I do not know how to explain this difference. 

{\bf Remark.} The ``type'' of $W$ in (6.2.7) depends not only on numbers $j,k$, but on a number $\mu$ from the equality 
$$pW \cap D_{g-1}= (\n Z/p)^{\mu}\eqno{(6.2.12)}$$
(see Appendix 2). I do not know what is the influence of $\mu$ on the above formulas. 
\medskip
\medskip
{\bf 6.3. Application to irreducible components of $T_{p,i}(V)$. }
\medskip
Apparently, the situation is similar to the one of section 4.3. We use notations of this section. Formulas 4.3.0, 4.3.1 hold for $\goth T_p=T_{p,i}$, (4.3.2) is rewritten as 
$$S_{g,i}= \bigcup _{l\in L} S'_{g,i}(l)\eqno{(6.3.2)}$$
and the conjecture (4.3.4) is rewritten as 
$$\tilde \goth z_{l} = \sum _{j,k} c_{jk\ l} \Phi_j\Phi_k (\tilde y_1)\eqno{(6.3.3)}$$
An analog of ${\pi'}^*_{j,l}$ is the following. For any $l \in L$, for any $j,k=0, \dots, g-1$ we consider the set $S^*_{g,i}(j,k) \cap S'_{g,i}(l)$ and the restriction of $\pi^*_{i;j,k}: S^{*}_{g,i}(j,k) \to G(j,k, g-1)$ on it. We denote this restriction by ${\pi'}^*_{i;j,k;l}: S^*_{g,i}(j,k) \cap S'_{g,i}(l) \to G(j,k, g-1)$. Analog of 4.3.5 is the following
\medskip
{\bf Condition 6.3.4.} For any $t\in G(j,k,g-1)$ the quantity of points in the fiber $({\pi'}^*_{i;j,k;l})^{-1}(t)$ is the same (does not depend on $t$). 
\medskip
We denote this quantity by $n_{ijkl}$. 
\medskip
{\bf Conjecture 6.3.5.} We have a formula
$$\tilde \goth z_l=\sum_{j,k=0}^{g-1} \frac{n_{ijkl}}{\alpha_{l}} z_{j,k}$$ 
Partition (6.3.2) for $g=3$, $\goth T_p=T_{p,1}$ is not known completely. We know $S'_{3,1}(l)$ for $l\in L_{\good}\times \Gal(K^p/K)$ and $l=j_0$ (notations of section 5). Since $g=3$, we have $D_{g-1}=D_2$ is 1-dimensional over $O(K_p)/p^2O(K_p)$, hence the phenomenon of section 4.5 does not hold for this case and we can expect that condition 6.3.4 is true for $g=3$, $\goth T_p=T_{p,1}$. 
\medskip
\medskip
{\bf Appendix 1. Calculation of $\tilde B$. }
\medskip
(2.2) shows that we can identify the basis $\goth B=\{e_1$, $e_2$, $e_3$, $e_4\}$ of $\tilde E_M$ with the basis $\{\goth f_{\emptyset}$, $\goth f_2$, $\goth f_3$, $\goth f_{23}\}$. In our case $\goth d=2$, i.e. $e_{\goth d+1}=e_3$. According formulas of Section 2 ((2.10), (2.27), (2.28) and others), we have: 
\medskip
$(\tilde y_1)_{l^{\infty}}$ --- the image of 
$\im (\tilde y_1)$ in $\tilde E_M=\tilde E(\n F_{p^2})_{l^{\infty}}$ is $\gamma e_3=\gamma \goth f_3$
\medskip
where $\gamma\in(\n Z/M)^*$. Application of (4.4.6) to the case $I=(3)$ gives us the following table of eigenvalues of elements of $\n H(\goth T)$ acting on eigenvector $e_3$ in $\tilde E_M$ (the second and fifth lines of the table are obtained by application of (4.4.6), the third and the sixth lines are obtained by application of (4.4.7): 

\medskip
\settabs 6 \columns
\+ Element of $\n H(\goth T)$ && $U_1U_2U_3$ & $U_1U_2V_3$ & $U_1V_2U_3$ & $U_1V_2V_3$ \cr 
\medskip
\+ Its eigenvalue && $a_0b_1b_2$ & $a_0b_1b_2b_3$ & $a_0b_1$ & $a_0b_1b_3$ \cr
\medskip
\+ Its residue modulo $M$ && $-1$ & $1$ & $-1$ & $1$ & \cr 
\medskip
\+ Element of $\n H(\goth T)$ && $V_1U_2U_3$ & $V_1U_2V_3$ & $V_1V_2U_3$ & $V_1V_2V_3$ \cr
\medskip
\+ Its eigenvalue && $a_0b_2$ & $a_0b_2b_3$ & $a_0$ & $a_0b_3$  \cr
\medskip
\+ Its residue modulo $M$ && $1$ & $-1$ & $1$ & $-1$ & \cr 
\medskip
Using (4.4.5) we get that the eigenvalues of $\Phi_i$, $i=0,1,2,3$, on $\tilde y_1$ are $$-1, \ \ 1, \ \ 1, \ \ -1\eqno{(A1.1)}$$ respectively (all formulas are in $\tilde E_M$). 

(4.2.2) for $g=3$ is the following: 

$$\matrix \Phi_0(\tilde y_1) & = & p^6z_0 \\ \Phi_1(\tilde y_1) & = & p^5z_0 + p^3z_1 \\ \Phi_2(\tilde y_1) & = & p^2z_1 + pz_2 \\ \Phi_3(\tilde y_1) & = & z_2 \endmatrix\eqno{(A1.2)}$$
(A1.1), (A1.2) and $p \equiv -1 \mod M$ imply that $z_0=z_2=-(\tilde y_1)_{l^{\infty}}$, $z_1=0$. Substituting these values in (4.3.10), (4.3.11) we get that $\kappa_p\equiv 0 \mod M$ and both $(\tilde y_p)_{l^{\infty}}$, $(\tilde y_{p,\bad})_{l^{\infty}}$ are 0. 
\medskip
\medskip
{\bf Appendix 2. Calculation of the quantity of poins in the fiber of $\pi^{*}_{1,i;j,k}$.}
\medskip
We use notations: $R=\n Z/p^2$, $B=(A_t)_{p^2}=R^6$. Let $e_1, \dots, e_6$ be an $R$-basis of $B$ such that the matrix of the skew form on $B$ in $e_1, \dots, e_6$ is $J_6$. Further, let $D=D_2 = < e_1, e_2>$, $W$ is an isotropic $R$-submodule of $B$ which is isomorphic to $R^2 \oplus (\n F_p)^2$ as an abstract module. We denote $W_4=W\cap pB$, $W_2=W_4^{\bot}=pW$. Finally, we denote $k-j$ by $\nu$. 
\medskip
We have: $W \cap D = R^j \oplus (\n F_p)^{\nu}$, $W_2 \cap D = \n F_p^{\mu}$, $j,{\nu},{\mu}$ are invariants of a given $W$ ($\mu$ of 6.2.12). 
\medskip
Problem: for each possible triples $j,{\nu},{\mu}$ find the quantity of $W$ with these invariants. 
\medskip
{\bf Lemma A2.1.} For a fixed $W_4$ there are $p^3$ $W$ such that its $W_4$ is the fixed one. 
\medskip
{\bf Proof.} I think that this quantity does not depend on a choice of $W_4$. If we take $W_4=e_1, \dots, e_4$ then it is possible always to choose 
$W=<e_2+p\alpha_{25}e_5+p\alpha_{26}e_6, e_3+p\alpha_{35}e_5+p\alpha_{36}e_6, pe_1, pe_4>$ where $(\alpha_{**})$ is a symmetric matrix, and this representation is unique. $\square$
\medskip
We consider in the following lemma the case of spaces over $\n F_p$. Let $B=\n F_p^6$, $D\subset B$ a fixed isotropic subspace of dimension 2, $W_2$ a variable isotropic subspace of dimension 2, and $W_4=W_2^{\bot}$. We denote $j_4= \dim W_4 \cap D$, $j_2= \dim W_2 \cap D$. 
\medskip
{\bf Lemma A2.2.} The quantity of $W_4$ with given $j_4, j_2$ is given by the following table: 

$$\matrix &&&&j_4&\\
&&&0&1&2\\
&&&&& \\
&0&&p^7&p^6 + 2p^5 + p^4 &0\\
j_2&1&&0&p^4+p^3&p^3 + 2p^2 + p\\
&2&&0&0&1 \endmatrix\eqno{(A2.3)}$$
\medskip
{\bf Proof.} Always we consider the quantity of $W_2$ with a base $x_1$, $x_2$. A fixed $W_2$ has $(p^2-1)(p^2-p)$ such bases. 
\medskip
(a) $j_4=0$. $W_4 \cap D =0 \iff W_2 \oplus D^{\bot}=B \iff W_2 \cap D^{\bot}=0$. There are $p^6-p^4$ possibilities for $x_1$; we have: $x_2 \in x_1^{\bot} - < x_1, D^{\bot}>$. It is easy to check that always $x_1^{\bot} \ne < x_1, D^{\bot}>$, i.e. there are always $p^5-p^4$ possibilities for $x_2$. 
\medskip
(b) $j_2=0$. There are $p^6-p^2$ possibilities for $x_1$; we have: $x_2 \in x_1^{\bot} - < x_1, D>$. 
\medskip
There are $2$ possibilities: 
\medskip
(1) $x_1^{\bot} \supset < x_1, D>$; 
\medskip
(2) $x_1^{\bot} \not\supset < x_1, D>$
\medskip
(1) $\iff x_1^{\bot} \supset D \iff x_1 \in D^{\bot}$. There are $p^4-p^2$ of such $x_1$, so there are $p^6-p^4$ of $x_1$ of type (2). The quantity of $x_2$ for $x_1$ of type (1) is $p^5-p^3$ and the quantity of $x_2$ for $x_1$ of type (2) is $p^5-p^2$. The desired quantity is $\frac{(p^4-p^2)(p^5-p^3)+(p^6-p^4)(p^5-p^2)}{(p^2-1)(p^2-p)}=p^7+p^6 + 2p^5 + p^4$. 
\medskip
(c) $j_4=2$. $W_4 \supset D \iff W_2 \subset D^{\bot}$. There are $p^4-1$ possibilities for $x_1$; we have: $x_2 \in (x_1^{\bot} \cap D^{\bot}) - < x_1>$. 
\medskip
There are $2$ possibilities: 
\medskip
(1) $x_1^{\bot} \supset D^{\bot}$; 
\medskip
(2) $x_1^{\bot} \not\supset D^{\bot}$
\medskip
(1) $\iff x_1 \in D$. There are $p^2-1$ of such $x_1$, so there are $p^4-p^2$ of $x_1$ of type (2). The quantity of $x_2$ for $x_1$ of type (1) is $p^4-p$ and the quantity of $x_2$ for $x_1$ of type (2) is $p^3-p$. The desired quantity is $\frac{(p^2-1)(p^4-p)+(p^4-p^2)(p^3-p)}{(p^2-1)(p^2-p)}=p^3+2p^2 + p+1$. 
\medskip
(d) $j_2=1$. Firstly we consider only $W_2$ such that $W_2 \cap D=<e_2>$, and we take $x_1=e_2$. So, $x_2=\alpha_1 e_1 + \alpha_2 e_2 + \alpha_3 e_3 + \alpha_4 e_4 + \alpha_6 e_6$. Condition on $x_2$: $(\alpha_3, \alpha_4, \alpha_6)\ne (0,0,0)$. There are $p^5-p^2$ of such $x_2$. For any $W_2$ there are $p^2-p$ of $x_2$-s that give us this $W_2$. So, the quantity of $W_2$ such that $W_2 \cap D = e_2$ is $\frac{p^5-p^2}{p^2-p}$ and the quantity of $W_2$ such that $j_2=1$ is $(p+1)\frac{p^5-p^2}{p^2-p}=p^4+2p^3 + 2p^2 + p$. 
\medskip
This gives us all entries of the table. $\square$ 
\medskip
Now for eash type of the table we calculate quantities of $W$ of different types. 
\medskip
(a). Case $j_2=2$. There are $p^3$ $W$ over it, they have form $W=<e_1+p\alpha_{14}e_4+p\alpha_{15}e_5, e_2+p\alpha_{24}e_4+p\alpha_{25}e_5, pe_3, pe_6>$ where $(\alpha_{**})$ is a symmetric matrix of lemma A2.1. We have: 
\medskip
$ W \cap D = (\n Z/p^2)^r \oplus {\n Z/p}^{2-r}$ where $r$ is the corank of $(\alpha_{**})$. 
\medskip
(b) Case $j_2=1, j_4=2$. I think that the quantities of $W$ over a $W_2$ with a given type of intersections with $D$ does not depend on a choose of $W_2$. For the case $W_2=p<e_2, e_3>$ a $W$ over it has a form $W=<e_2+p\alpha_{25}e_5+p\alpha_{26}e_6, e_3+p\alpha_{35}e_5+p\alpha_{36}e_6, pe_1, pe_4>$ where $(\alpha_{**})$ is a symmetric matrix of lemma A2.1. We have: 

$$\alpha_{25}=\alpha_{26}=0 \iff W \cap D = \n Z/p^2 \oplus \n Z/p$$
there are $p$ such spaces $W$, and 
$$(\alpha_{25},\alpha_{26})\ne(0,0) \iff W \cap D = (\n Z/p)^2$$
there are $p^3-p$ such spaces $W$.
\medskip
(c) Case $j_2=1, j_4=1$. For the case $W_2=p<e_2, e_4>$ a $W$ over it has a form $W=<e_2+p\alpha_{21}e_1+p\alpha_{25}e_5, e_4+p\alpha_{41}e_1+p\alpha_{45}e_5, pe_3, pe_6>$ where $(\alpha_{**})$ is a symmetric matrix of lemma A2.1. We have: 

$$\alpha_{25}=0 \iff W \cap D = \n Z/p^2$$
there are $p^2$ such spaces $W$, and 
$$\alpha_{25}\ne 0 \iff W \cap D = \n Z/p$$
there are $p^3-p^2$ such spaces $W$. $\square$
\medskip
So, we have the following 
\medskip
{\bf Theorem A2.4.} The quantities of $W$ with invariants $(i,{\nu},{\mu})$ are the following: 
\medskip
\settabs 8 \columns
\+ $(j,{\nu},{\mu})$ & $W \cap D $ && $j_2$ & $j_4$ & quantity && \cr
\medskip
\+ 0,0,0 & 0&& 0 & 0 & $p^{10}$ &&\cr
\medskip
\+ 0,1,0 & $\n Z/p$ && 0 & 1& $p^9+2p^8+p^7$&&\cr
\medskip
\+ 0,1,1 & $\n Z/p$ && $1$ & 1& $p^7-p^5$&&\cr
\medskip
\+ 1,0,1 & $\n Z/p^2$ && $1$ & 1&$p^6+p^5$&&\cr
\medskip
\+ 0,2,0 & $\n Z/p \oplus  \n Z/p$ && 0 & 2& 0 \cr
\medskip
\+ 0,2,1 & $\n Z/p \oplus  \n Z/p$ && $1$ & 2& $(p-1)p^2(p+1)^3$\cr
\medskip
\+ 0,2,2 & $\n Z/p \oplus  \n Z/p$ && $2$ & 2&$p^3-p^2$\cr
\medskip
\+ 1,1,1 & $\n Z/p^2 \oplus  \n Z/p$ && $1$ & 2&$p^4+2p^3+p^2$\cr
\medskip
\+ 1,1,2 & $\n Z/p^2 \oplus  \n Z/p$ && $2$ & 2&$p^2-1$\cr
\medskip
\+ 2,0,2 & $\n Z/p^2 \oplus  \n Z/p^2$ && $2$ & 2&1 \cr
\medskip
\medskip
{\bf Appendix 3. Idea of finding of $\goth x_p$ (based on [Z85]).} \nopagebreak
\medskip
Recall that a 1-dimensional analog of the problem of finding of $\goth x_p$ is a calculation of Gross-Kohnen-Zagier ([Z85], [GKZ87]) of ratios of Heegner points corresponding to different imaginary quadratic fields on a fixed elliptic curve over $\n Q$. We reformulate here the crucial theorem of [Z85] for a general situation, and we consider possibilities of application of this theorem. 
\medskip
Let $X$ be the variety under consideration. It is defined over $\n Q$. 
\medskip
Let $C^r(X)$ be the group of cycles on $X$ of codimension $r$, 
\medskip
$Ch: C^r(X) \to H^{2r}(X_{\n Q}, \n Z_l(r))$ the cycle map, 
\medskip
$j_0: H^{2r}(X_{\n Q}, \n Z_l(r)) \to H^0(\Gal (\bar \n Q/ \n Q), H^{2r}(X_{\bar \n Q}, \n Z_l(r)))$ and 
\medskip
$j_1: \Ker j_0 \to H^1(\Gal (\bar \n Q/ \n Q), H^{2r-1}(X_{\bar \n Q}, \n Z_l(r)))$ maps coming from the spectral sequence, 
\medskip
$\cl = j_0 \circ Ch: C^r(X) \to H^0(\Gal (\bar \n Q/ \n Q), H^{2r}(X_{\bar \n Q}, \n Z_l(r)))$ the cycle map, 
\medskip
$C^r(X)_0$ its kernel and 
\medskip
$\cl_1 = j_1 \circ Ch: C^r(X)_0 \to H^1(\Gal (\bar \n Q/ \n Q), H^{2r-1}(X_{\bar \n Q}, \n Z_l(r)))$ the Abel-Jacobi map. 
\medskip
We denote its kernel by $C^r(X)_1$. 
\medskip
We need a criterion that a cycle $C \in C^r(X)$ belongs to $C^r(X)_1$. 
\medskip
{\bf Theorem.} ([Z85]). Let there exist a variety $Y$, an inclusion $i: X \hookrightarrow Y$ and a cycle $Z \in C^r(Y)$ such that: 
\medskip
(a) $H^{2r-1}(Y_{\bar \n Q}, \n Z_l(r)))=0$; 
\medskip
(b) $i^*(Z)=C$; 
\medskip
(c) $Z \in C^r(Y)_0$. 
\medskip
Then $C \in C^r(X)_1$. $\square$
\medskip
{\bf Problem.} Apply this theorem to the case $X=X_0(N)$ is a smooth completion of the Siegel 6-fold of level $N$, $C$ is an appropriate linear combination of $V$, $V_{p,i}$, $V_{p,j}$ from (1.6) for the case $\goth T_p=T_{p,1}$, and cycles with support at infinity of $X$. ($N$ is denoted by $M$ in [Z85]). 
\medskip
Ideas of answer. 
\medskip
We can take: 
\medskip
(a) $Y$ is a Hilbert-Siegel variety ($g=2$, $m=6$, notations of [Sh63]); 
\medskip
(b) $Y=X_0(1) \times X_0(1)$, where $X_0(1)$ is a smooth completion of the Siegel 6-fold of level 1. 
\medskip
For both $Y$ there exists an inclusion $i: X \hookrightarrow Y$ ($N$ and the real quadratic field in the case (a) must satisfy a certain condition). 

If we were able to find $Z_0 \in C^4(Y)_0$ such that (b) of the lemma holds then we could expect that $Z_{p,1}$ --- the union of irreducible components of $T_{p,?}(Z_0)$ which are defined over $K$ (here $T_{p,?}$ is a $p$-Hecke correspondence on $Y$) --- can be the desired $Z$. 

Maybe we try to find a Hodge class (i.e. an element of $H^4(Y, \n Z) \cap H^{2,2}(Y)$) such that its intersection with the class of $X$ gives us the class of $C$? 
\medskip
{\bf References. }
\medskip
[BK90] Bloch S., Kato L., L-functions and Tamagava numbers of 
motives. In: The Grothendieck Festschrift, v.1, p. 333 - 400. Boston 
Basel Stuttgart: Birkhauser, 1990

[D71] Deligne P. Travaux de Shimura. Lect. Notes in Math., 1971, 
v.244, p. 123 - 165. Seminaire Bourbaki 1970/71, Expos\'e 389. 

[G91] Gross B.H., Kolyvagin's work on modular elliptic 
curves. In: $L$-functions and Arithmetic. Proceedings of the Durham 
Symposium, July, 1989. Cambridge Univ. Press, 1991. p. 235 - 256. 

[GZ86] Gross B.H., Zagier D.B., Heegner points and derivatives of $L$-
series. Invent. Math., 1986, v.84, N.2, p. 225 - 320.

[GKZ87] Gross B.H., Kohnen W., Zagier D.B., Heegner points and
derivatives 
of $L$-series, II.  Math. Ann., 1987, v.278.

[K89] Kolyvagin V.A., Finiteness of $E(\n Q )$ and Sh$(E,\n Q )$ 
for a subclass of  Weil curves. Math. USSR Izvestiya, 1989, v. 32, 
No. 3, p. 523 - 541 

[K90] Kolyvagin V.A. Euler systems. In: The Grothendieck 
Festschrift, v.2. Boston Basel Stuttgart: Birkhauser 1990, p. 435 - 483   

[KL92] Kolyvagin V.A., Logachev D.Yu. Finiteness of Sh over totally 
real fields. Math. USSR Izvestiya, 1992, v. 39, No. 1, p. 829 - 853

[L01] Logachev D., Action of Hecke correspondences on Heegner curves on a Siegel threefold. J. of Algebra, 2001, v. 236, N. 1, p. 307 - 348 

[L04] Logachev D., Relations between conjectural eigenvalues of Hecke operators on submotives of Siegel varieties. http://arxiv.org/ps/math.AG/0405442
 
[L05] Logachev D., Action of Hecke correspondences on subvarieties of Shimura varieties. http://arxiv.org/ps/math.AG/0508204

[N92] Nekov\'a\v r J. Kolyvagin's method for Chow groups of 
Kuga-Sato varieties. Inv. Math. 1992, v. 107, p. 99 - 125. 

[Sh63] Shimura G., On analytic families of polarized abelian 
varieties and automorphic functions. Annals of Math., 1963, v. 78, No. 1, p. 
149 - 192

[Z85] Zagier, Don,  Modular points, modular curves, modular surfaces and modular forms. Workshop Bonn 1984 (Bonn, 1984), 225 - 248, Lecture Notes in Math., 1111, Springer, Berlin, 1985.
\medskip
E-mail: logachev\@ usb.ve
\enddocument